\documentclass[review,12pt]{elsarticle}

\usepackage[cmex10]{amsmath}
\usepackage{url}
\usepackage{amssymb}
\usepackage{mathrsfs}
\usepackage{comment}
\usepackage{graphicx}
\usepackage{rotating}
\usepackage{bm}
\usepackage{array}
\usepackage{tabularx,ragged2e,booktabs,caption}

\newcommand{\iffm}{\Longleftrightarrow}
\newcommand{\bs}[1]{\boldsymbol{#1}}

\newcommand{\gr}[1]{\mathbf{#1}}

\newcommand{\E}{\mathbb{E}}

\newcommand{\En}{\mathrm{E}}
\newcommand{\Jn}{\mathrm{J}}

\newcommand{\wW}{\widetilde{W}}

\newcommand{\Rx}{R_{\x}}
\newcommand{\wRx}{R_{\w{\x}}}

\newcommand{\Ry}{R_{\y}}
\newcommand{\wwRy}{\w{R}_{\y}}
\newcommand{\Rn}{R_{\n}}
\newcommand{\wRn}{R_{\w{\n}}}
\newcommand{\Rxy}{R_{\x\y}}
\newcommand{\wwwRxy}{\w{R}_{\wx\y}}
\newcommand{\Ryx}{R_{\y\x}}
\newcommand{\wwwRyx}{\w{R}_{\y\wx}}
\newcommand{\Rxn}{R_{\x\n}}

\newcommand{\Pre}{\gr{P}}

\newcommand{\n}{\gr{n}}
\newcommand{\wn}{\w{\gr{n}}}

\newcommand{\x}{\gr{x}}
\newcommand{\wx}{\w{\gr{x}}}
\newcommand{\z}{\gr{z}}

\newcommand{\y}{\gr{y}}

\newcommand{\w}[1]{\widetilde{#1}}

\newcommand{\sm}{\mathbb{R}^m}

\newcommand{\nn}{\mathbb{R}^{n\times n}}

\newcommand{\mm}{\mathbb{R}^{m\times m}}
\newcommand{\mr}{\mathbb{R}^{m\times r}}

\newcommand{\mn}{\mathbb{R}^{m\times n}}
\newcommand{\nm}{\mathbb{R}^{n\times m}}

\newcommand{\tr}{\mathbb{R}^{t\times r}}

\newcommand{\zt}{\textrm}

\newcommand{\ra}[1]{\mathcal{R}\left({#1}\right)}

\newcommand{\p}{\parallel}

\newcommand{\si}{\sigma}
\newcommand{\be}{\beta}
\newcommand{\ep}{\epsilon}
\newcommand{\ga}{\gamma}
\newcommand{\Si}{\Sigma}

\newcommand{\Ga}{\Gamma}

\newcommand{\up}{\upsilon}
\newcommand{\De}{\Delta}

\newcommand{\de}{\delta}

\newcommand{\no}{\noindent}
\newcommand{\pd}{\no\hspace{2em}{\itshape Proof: }} 
\newcommand{\ds}{\displaystyle}

\newtheorem{theorem}{Theorem}

\newtheorem{lemma}{Lemma}
\newtheorem{corollary}{Corollary}
\newtheorem{proposition}{Proposition}

\newtheorem{definition}{Definition}
\newtheorem{remark}{Remark}
\newtheorem{fact}{Fact}
\newcolumntype{C}[1]{>{\Centering}m{#1}}

\journal{Journal of The Franklin Institute}

\begin{document}

\begin{frontmatter}



  \title{Reduced-Rank Estimation for Ill-Conditioned Stochastic Linear Model with High Signal-to-Noise Ratio}
  \author[TP,TP2]{Tomasz Piotrowski\corref{cor1}}
  \ead{tpiotrowski@is.umk.pl}
  \cortext[cor1]{Corresponding author.}
  \fntext[1]{Published: \url{https://doi.org/10.1016/j.jfranklin.2016.05.007}}
  \address[TP]{Department of Informatics,\\
    Faculty of Physics, Astronomy and Informatics,\\
    Nicolaus Copernicus University,\\
    Grudziadzka 5, 87-100 Torun, Poland}
  \address[TP2]{Interdisciplinary Center of Modern Technologies,\\
    Nicolaus Copernicus University,\\
    Wilenska 4, 87-100 Torun, Poland}

  \author[IY]{Isao Yamada}
  \ead{isao@sp.ce.titech.ac.jp}
  \address[IY]{Department of Communications and Computer Engineering,\\
    Tokyo Institute of Technology,\\
    Tokyo 152-8550, Japan}


  \begin{abstract}
    Reduced-rank approach has been used for decades in robust linear estimation of both deterministic and random vector of parameters in linear model $\y=H\x+\sqrt{\ep}\n.$ In practical settings, estimation is frequently performed under incomplete or inexact model knowledge, which in the stochastic case significantly increases mean-square-error (MSE) of an estimate obtained by the linear minimum mean-square-error (MMSE) estimator, which is MSE-optimal among linear estimators in the theoretical case of perfect model knowledge. However, the improved performance of reduced-rank estimators over MMSE estimator in estimation under incomplete or inexact model knowledge has been established to date only by means of numerical simulations and arguments indicating that the reduced-rank approach \emph{may} provide improved performance over MMSE estimator in certain settings. In this paper we focus on the high signal-to-noise ratio (SNR) case, which has not been previously considered as a natural area of application of reduced-rank estimators. We first show explicit sufficient conditions under which familiar reduced-rank MMSE and truncated SVD estimators achieve lower MSE than MMSE estimator if singular values of array response matrix~$H$ are perturbed. We then extend these results to the case of a generic perturbation of array response matrix~$H$, and demonstrate why MMSE estimator frequently attains higher MSE than reduced-rank MMSE and truncated SVD estimators if~$H$ is ill-conditioned. The main results of this paper are verified in numerical simulations.
  \end{abstract}

  \begin{keyword}
    robust linear estimation \sep reduced-rank estimation \sep stochastic MV-PURE estimator \sep array signal processing
  \end{keyword}

\end{frontmatter}


\section{Introduction}
Linear estimation of random vector of parameters $\x$ in a stochastic linear model $\y=H\x+\sqrt{\ep}\n$ under mean-square-error (MSE) criterion continues to be one of the central problems of signal processing \cite{Luenberger1969,Kailath2000}, e.g., brain signal processing \cite{Cichocki2002,Piotrowski2013}, wireless communications \cite{Wang2004}, and array signal processing \cite{VanTrees2002}, \cite[Sec.1.2]{Pezeshki2010}. The linear minimum mean-square-error estimator (MMSE, often called the Wiener filter) \cite{Luenberger1969,Kailath2000} achieves the lowest MSE among linear estimators if covariance matrix $\Ry$ of $\y$ and cross-covariance matrix $\Rxy$ between $\x$ and $\y$ are known exactly. However, its performance degrades significantly if only their estimates are available, which happens for example if $H$ deviates from the one assumed, see, e.g., \cite{Eldar2004, Mittelman2010, Kalantarova2012, Zachariah2014} and references therein.

To alleviate this problem, a myriad of solutions have been proposed over the years, with the reduced-rank approach being one of the most promising approaches demonstrating much improved robustness to imperfect model knowledge compared to theoretically MSE-optimal MMSE estimator \cite{Scharf1991,Yamashita1996,Stoica1996,Hua2001}.\footnote{The reduced-rank approach has also found use in robust deterministic least-squares estimation \cite{Marquardt1970,Chipman1999,Werner2006}, alongside methods such as total least squares \cite{Golub1980,VanHuffel1991,Sayed1998}.} However, this improved performance of reduced-rank estimators has only been demonstrated via numerical simulations or arguments indicating that the reduced-rank approach \emph{may} provide improved performance over MMSE estimator in certain settings. Such situation renders reduced-rank estimation confined to specific applications, where it is expected to provide additional robustness to imperfect model knowledge compared to MMSE estimator. In particular, the  reduced-rank estimators are usually not applied in high signal-to-noise ratio (SNR) settings, and the body of research in this area is very limited.

We fill this gap by analyzing in detail the MSE performance of the familiar reduced-rank MMSE and truncated SVD estimators for high SNR case and realistic small perturbations of an ill-conditioned array response matrix. Namely, under the assumption of spatially white random vector of parameters $\x$ and spatially white noise $\n$, we consider the reduced-rank MMSE estimator \cite{Scharf1991,Hua2001} and the truncated SVD estimator, which in such settings is equivalent to the stochastic MV-PURE estimator \cite{Piotrowski2009a,Piotrowski2014} (see \cite{Yamada2006,Piotrowski2008,Piotrowski2009a,Piotrowski2009b,Piotrowski2012b,Piotrowski2013,Yamagishi2013,Piotrowski2014} for details on the MV-PURE framework). We provide explicit sufficient conditions for reduced-rank MMSE and truncated SVD estimators to achieve lower MSE than MMSE estimator if singular values of an ill-conditioned array response matrix $H$ are mildly perturbed and the SNR is high. We also compare the performances of the reduced-rank MMSE and truncated SVD estimators in such settings. We then extend the analysis to generic perturbation of $H$ and show that the MSEs of reduced-rank MMSE and truncated SVD estimators are approximately equal to the respective MSEs obtained for the simplified perturbation model under mild perturbation of the array response matrix. Moreover, we show that in such settings the MSE of MMSE estimator is likely to be larger than the MSE of the reduced-rank MMSE estimator. A case study for square Gaussian distributed $H$ is also given, illustrating the applicability of the derived results. The main results of this paper are verified in numerical simulations.

The paper is organized as follows. In Section \ref{preeli} we introduce necessary preliminaries. In Section \ref{defining} we establish key definitions regarding ill-conditioned stochastic linear model with high signal-to-noise ratio in our sense. In Section \ref{simple} we derive explicit sufficient conditions for reduced-rank MMSE and truncated SVD estimators to achieve lower MSE than MMSE estimator if singular values of array response matrix are mildly perturbed and SNR is high. In Section~\ref{ext} we extend the results of Section \ref{simple} to the case of generic perturbations of array response matrix $H.$ In Section \ref{ne} we provide a numerical example illustrating the main results of the paper. Section \ref{conclusion} concludes our work.

For readers' convenience, the necessary mathematical facts used in this paper are stated in \ref{kru}. 

A short preliminary version of the paper was presented at conference \cite{Piotrowski2012}.

\section{Preliminaries} \label{preeli}
\subsection{Stochastic Linear Model}
Consider the stochastic linear model of the form:
\begin{equation} \label{general}
  \y=H\x+\sqrt{\ep}\n,
\end{equation}
where $\y, \x, \n$ are random vectors representing observed signal, signal to be estimated, and additive noise, respectively, $H\in\nm$ is a known matrix of rank $m$, and $\ep>0$ is a known constant representing noise power. We assume that $\x$ and $\n$ have zero mean, are uncorrelated: $\Rxn=0\in\mn$, and white: 
\begin{equation} \label{biel}
  \Rx=I_m\zt{ and }\Rn=I_n,
\end{equation}
where by $\Rxn$ we denote the cross-covariance matrix of $\x$ and $\n$, by $\Rx$ and $\Rn$ we denote the covariance matrices of $\x$ and $\n$, respectively, and by $I_s$ the identity matrix of size $s.$ Note that from our assumptions $\Ry=(HH^t+\ep I_n)\succ 0$ and $\Ryx=H$ are available, where by $X\succ 0$ we mean that a square matrix $X$ is positive-definite.

We denote singular value decomposition (SVD) of $H$ by
\begin{equation} \label{svdtruH}
  H=M\Ga N^t,
\end{equation}
where $M=(m_1,\dots,m_n)$, $N=(n_1,\dots,n_m)$, with distinct singular values $\ga_i,i=1,\dots,m$, organized in decreasing order.

The norm of a random vector $\z$ is defined as $\p\z\p=\sqrt{tr[\E(\z\z^t)]}$, see, e.g., \cite{Luenberger1969}. Based on this, we define signal-to-noise ratio (SNR) as:
\begin{equation} \label{SNR}
  SNR={{\p H\x\p^2}\over{\p\sqrt{\ep}\n\p^2}}={{tr\left[\E\left[H\x\x^tH^t\right]\right]}\over{tr\left[\E\left[\sqrt{\ep}\n\sqrt{\ep}\n^t\right]\right]}}={{tr[HH^t]}\over{n\ep}}=(n\ep)^{-1}\sum_{i=1}^m\ga_i^2.
\end{equation}

The stochastic linear model (\ref{general}) has found a widespread use in signal processing \cite{Luenberger1969,Kailath2000}. In certain applications, such as in wireless communications, it is often natural to assume that the signal $\x$ to be estimated is white or can be whitened a priori. Then, if the noise covariance matrix is known, the assumptions (\ref{biel}) will be satisfied by premultiplying the observed signal with $\Rn^{-1/2}$, which is the inverse of the positive-definite square root of $\Rn.$\footnote{For every natural $k\geq 1$ and every positive definite matrix $X\succ 0$ there exists a unique positive definite matrix $Y\succ 0$ such that $Y^k=X$ \cite[Th.7.2.6, p.405]{Horn1985}.} Moreover, the recent work \cite{Mizoguchi2014} showed that any hypercomplex linear system, including complex widely linear model, is algebraically equivalent to model (\ref{general}). In such a case, the assumptions (\ref{biel}) will be satisfied, e.g., if the hypercomplex-valued desired signal is designed such that its real-valued representation satisfies $\Rx=I_m$, and the noise is white circular Gaussian.   

On the other hand, the exact knowledge on model (\ref{general}) is usually unavailable in practice. We explore consequences of this fact in this paper, and consider the situation where one has only imperfect value $H+\Delta H$ in place of the true value $H$ of the array response matrix in (\ref{general}) and therefore one has to use the following misguided linear model [\emph{cf.} (\ref{general})]:
\begin{equation} \label{bubu}
  \y=(H+\De H)\wx+\sqrt{\ep}\wn,
\end{equation}
where $\De H\in\nm$ represents the perturbation of the array response matrix $H$ such that $H+\De H$ is also of rank~$m.$ As a consequence of replacing $H$ with $H+\De H$ in the model definition, we note that model (\ref{bubu}) employs random variables $\wx$ and $\wn$ to explain observation $\y$, which differ from $\x$ and $\y$ used in model (\ref{general}) to explain \emph{the same} observation $\y.$ Clearly, we may still assume that $\w{\x}$ and $\w{\n}$ are uncorrelated, have zero mean and satisfy
\begin{equation} \label{biel_p}
  \wRx=I_m\zt{ and }\wRn=I_n.
\end{equation}
However, their realizations may be different from the corresponding $\x$ and $\n$ in model (\ref{general}), as both (\ref{general}) and (\ref{bubu}) yield the same observation $\y$, but use different array response matrices. In particular, the assumed form of the covariance matrix of $\y$ will also differ from the actual form. Namely, in such a situation, we often have to use the erroneous covariance matrix 
\begin{equation} \label{Ry_pe}
  \wwRy=(H+\De H)(H+\De H)^t+\ep I_n.
\end{equation}
Similarly, in place of $\Ryx=H$ one has access only to:
\begin{equation} \label{Ryx_pe}
  \wwwRyx=H+\De H.
\end{equation}

We denote singular value decomposition (SVD) of $H+\De H$ by
\begin{equation} \label{svdH}
  H+\De H=U\Si V^t,
\end{equation}
where $U=(u_1,\dots,u_n)$, $V=(v_1,\dots,v_m)$, with distinct singular values $\si_i,i=1,\dots,m$, organized in decreasing order, and by $V_r=(v_1,\dots,v_r)\in\mr$ we denote first $r$ columns of $V.$ 

\subsection{Linear Estimation - MMSE, $r$-MMSE and $r$-SVD Estimators} \label{linest}
For the linear estimators considered in this paper, we assume that they are constructed based on available perturbed model (\ref{bubu}), but they do not have access to the true array response matrix $H$ of the unperturbed model~(\ref{general}). Thus, the best possible effort to estimate $\x$ from $\y$ is to estimate instead $\wx$ from $\y$ in model (\ref{bubu}) as
\begin{equation} \label{estimate}
  \widehat{\wx}=\wW\y,
\end{equation}
where $\wW\in\mn$ is a matrix called here an estimator. The tilde over $\wW$ emphasizes that it is constructed based on available misguided model (\ref{bubu}).

However, we aim at estimating $\x$ in model (\ref{general}). Therefore, we should evaluate the mean-square-error of an estimate in (\ref{estimate}) in terms of model (\ref{general}). More precisely, the MSE of $\widehat{\wx}$ in~(\ref{estimate}) as an estimate of $\x$ from observation $\y$ given in (\ref{general}) as $\y=H\x+\sqrt{\ep}\n$ is of the form:
\begin{multline} \label{mse}
  \Jn(\wW)=tr\left[\En[(\widehat{\wx}-\x)(\widehat{\wx}-\x)^t]\right]=tr\left[\En[(\wW\y-\x)(\wW\y-\x)^t]\right]=\\
  tr\left[\En[\wW\y\y^t\wW^t]\right]-2tr\left[\En[\wW\y\x^t]\right]+tr\left[\En[\x\x^t]\right]=\\
  tr[\wW\Ry \wW^t]-2tr[\wW\Ryx]+tr[\Rx]=\\ tr[\wW(HH^t+\ep I_n) \wW^t]-2tr[\wW H]+m.
\end{multline}

We introduce now the estimators considered in this paper. The linear minimum mean-square-error estimator (MMSE) (often called the Wiener filter) designed based on model (\ref{bubu}) is of the following form \cite{Luenberger1969,Kailath2000,VanTrees2002}:
\begin{equation}
  \wW_{MMSE}=\wwwRxy\wwRy^{-1}.
\end{equation}
We may express the MMSE estimator in view of (\ref{Ry_pe}) and (\ref{Ryx_pe}) as:
\begin{equation} \label{mmseish}
  \wW_{MMSE}=(H+\De H)^t[(H+\De H)(H+\De H)^t+\ep I_n]^{-1}=
  V\Si^t(\Si\Si^t+\ep I_n)^{-1}U^t.
\end{equation}

For a given rank constraint $1\leq r<m$, the reduced-rank MMSE estimator \cite{Scharf1991,Hua2001} can be expressed in our settings as:
\begin{equation} \label{rrmmseish}
  \wW_{r-MMSE}=V_rV_r^t\wW_{MMSE}=VI_m^r\Si^t(\Si\Si^t+\ep I_n)^{-1}U^t,
\end{equation}
where $I_m^r\in\mm$ contains as its $r\times r$ principal submatrix the identity matrix of size $r$ and zeros elsewhere.

Similarly, for a given rank constraint $1\leq r<m$, the truncated SVD estimator, which in our settings is the simplest example of the stochastic MV-PURE estimator \cite{Piotrowski2009a,Piotrowski2014} (see \cite{Yamada2006,Piotrowski2008,Piotrowski2009a,Piotrowski2009b,Piotrowski2012b,Piotrowski2013,Yamagishi2013,Piotrowski2014} for details on the MV-PURE framework) is given by:
\begin{equation} \label{mvpform_s}
  \wW_{r-SVD}=V_rV_r^t(H+\De H)^\dagger=VI_m^rV^t(H+\De H)^\dagger=V\Si_r^\dagger U^t,
\end{equation}
where $(H+\De H)^\dagger$ is the Moore-Penrose pseudoinverse of $H+\De H$ \cite{BenIsrael2003}, and where $\Si_r^\dagger\in\mn$ is of the following form
\begin{equation}
  \Si_r^\dagger=\left(
  \begin{array}{ccccccc} 
    \si_1^{-1} & 0 & 0 & 0 & 0 & \dots & 0\\
    0 & \si_2^{-1} & 0 & 0 & 0 & \dots & 0\\
    0 & 0 & \ddots & 0 & 0 & \dots & 0\\
    0 & 0 & 0 & \si_r^{-1} & 0 & \dots & 0\\
    0 & 0 & 0 & 0 & 0 & \dots & 0\\
    \vdots & \vdots & \vdots & \vdots & \vdots & \vdots & \vdots\\
    0 & 0 & 0 & 0 & 0 & \dots & 0\\
  \end{array} \right)\in\mn.
\end{equation}.

\section{Ill-Conditioned Stochastic Linear Model with High Signal-to-Noise Ratio} \label{defining}
If we look at the definition of SNR in (\ref{SNR}), it is clear that in order to define an ill-conditioned stochastic linear model (\ref{general}) with high signal-to-noise ratio, we should clearly define the interplay between the noise power $\ep$ and the singular values $\ga_i$ of $H.$ Moreover, the \emph{ill-conditioned} part of the definition should reflect the fact that for a certain $r$ such that $1\leq r<m$, the trailing singular values $\{\ga_{r+1},\dots,\ga_m\}$ must be much smaller than the leading singular values $\{\ga_1,\dots,\ga_r\}.$ Similarly, the phrase \emph{high signal-to-noise ratio} implies that at least some of the the leading singular values $\{\ga_1,\dots,\ga_r\}$ must be substantially larger than the noise power $\ep$ to obtain a sufficiently high value of $SNR$~in~(\ref{SNR}) for all reasonable ratios $n/m\geq 1.$ In such a case, we may assume that the leading singular values $\{\ga_1,\dots,\ga_r\}$ are well separated from each other which occurs frequently in practice, see, e.g., the probability density function of the distribution of singular values in the case of Gaussian-distributed~$H$ in (\ref{pdf}) in Section \ref{cs}.

We should also expect that the perturbed model (\ref{bubu}) has approximately the same properties as the ones listed above. Therefore, it will be more natural to introduce a definition of a corresponding pair of models (\ref{general}) and (\ref{bubu}). Such approach will enable us not only to define the terms \emph{ill-conditioned stochastic linear model} and \emph{high signal-to-noise ratio} in our sense, but will also enable us to clearly define the interplay between model (\ref{general}) and its perturbed version (\ref{bubu}).

In order to formalize the above considerations, we prove first the following proposition which gives bounds on the maximum deviation of $\si_i$ from $\ga_i.$
\begin{proposition} \label{YamadaSensei}
  Let us consider SVD of $H$ of the form (\ref{svdtruH}) and SVD of $H+\De H$ of the form (\ref{svdH}), and let $i\in\{1,\dots,m\}.$ If:
  \begin{equation} \label{YSasu}
    \ga_i>\p\De H\p_2,
  \end{equation}
  then:
  \begin{equation}
    \si_i\in\left[\ga_i-\p\De H\p_2,\ \ga_i+\sqrt{2}\p\De H\p_2\right].
  \end{equation}
\end{proposition}
\pd See \ref{pd_YamadaSensei}.

In view of the discussion in the first two paragraphs of this section, for simplicity, we assume (\ref{YSasu}) is satisfied for $i=1,\dots,r$ for a certain $r$ such that $1\leq r<m.$ Based on this observation, we introduce the following definition.
\begin{definition} \label{master_def}
  Let $H$ have SVD of the form (\ref{svdtruH}), and $H+\De H$ have SVD of the form (\ref{svdH}). Moreover, let $r$ be such that $1\leq r<m.$ Then, we say that the singular values $\{\si_1,\dots,\si_r\}$ of $H+\De H$ and  $\{\ga_1,\dots,\ga_r\}$ of $H$ are pairwise separated if $\ga_r>\p\De H\p_2$ is satisfied and the intervals
  \begin{equation} \label{beauty2}
    \left[\ga_i-\p\De H\p_2,\ \ga_i+\sqrt{2}\p\De H\p_2\right] 
  \end{equation}
  are disjoint for $i=1,\dots,r.$ 
\end{definition}

\begin{remark}
  We note that for a given set of singular values $\ga_1,\dots,\ga_r$ of $H$, Definition~\ref{master_def} imposes an upper bound on the size of perturbation $\De H$ by insisting that assumption (\ref{YSasu}) is satisfied for $i=1,\dots,r$ and the intervals in (\ref{beauty2}) are disjoint for $i=1,\dots,r$, which implies $\si_1>\si_2>\dots>\si_r\geq\ga_r-\p\De H\p_2>0.$
\end{remark}

We may introduce now the following definition which captures the key points of the first paragraph of this section. More precisely, condition~\ref{c1}. of Definition \ref{highill} below ensures that singular values of $H$ and $H+\De H$ are pairwise separated, while conditions \ref{c2}. and \ref{cend}. introduce the distribution of singular values of $H$ and $H+\De H$ which yields ill-conditioned stochastic linear models (\ref{general}) and (\ref{bubu}) with high signal-to-noise ratio.
\begin{definition} \label{highill}
  Consider stochastic linear models (\ref{general}) and (\ref{bubu}). We say that (\ref{general}) and (\ref{bubu}) constitute a corresponding pair of ill-conditioned stochastic linear models with high signal-to-noise ratio of gap $\kappa>1$ under perturbation $\De H$ if for a certain $r$ such that $1\leq r<m$:
  \begin{enumerate}
  \item singular values $\{\ga_1,\dots,\ga_r\}$ of $H$ and $\{\si_1,\dots,\si_r\}$ of $H+\De H$ are pairwise separated, \label{c1}
  \item singular values $\{\ga_1,\dots,\ga_r\}$ of $H$ and $\{\si_1,\dots,\si_r\}$ of $H+\De H$ are greater than $\kappa\sqrt{\ep}$, \label{c2}
  \item singular values $\{\ga_{r+1},\dots,\ga_m\}$ of $H$ and $\{\si_{r+1},\dots,\si_m\}$ of $H+\De H$ are at most $\sqrt{\ep}.$ \label{cend}
  \end{enumerate}
\end{definition}

\begin{remark}
  In the subsequent sections we will consider models satisfying Definition~\ref{highill}. The reduced-rank estimators will be of rank $r$, where $r$ is given in the sense of Definition \ref{highill}.
\end{remark}

\section{Performance Gain for Perturbed Singular Values Model} \label{simple}
In this section we assume that SVD of the true array response matrix $H$ in~(\ref{general}) is of the following form:
\begin{equation} \label{asu}
  H=U\Ga V^t.
\end{equation}
Thus, we may set $M=U$ and $N=V$ in (\ref{svdtruH}), and evaluate the effects of perturbation of singular values of the array response matrix.\footnote{We note that the assumption of distinct singular values implies that if $\hat{M}\Ga\hat{N}^t$ is another SVD of $H$ [\emph{cf.} (\ref{svdtruH})] with $\hat{M}=(\hat{m}_1,\dots,\hat{m}_n)$, $\hat{N}=(\hat{n}_1,\dots,\hat{n}_m)$, then we must have either $\hat{m}_i=m_i$ and $\hat{n}_i=n_i$, or $\hat{m}_i=-m_i$ and $\hat{n}_i=-n_i$ for any $i\in\{1,\dots,m\}.$ This fact follows immediately from \cite[Th.3.1.1', p.147]{Horn1991} and the fact that for a real matrix we may only consider real singular value decompositions.} The analysis in this section will give a useful insight to deal with the case of an arbitrary
perturbation satisfying the conditions in Definition \ref{highill}.

\subsection{Direct Evaluation of Mean-Square-Error} \label{direct}
We begin with the following proposition.

\begin{proposition} \label{s_expr}
  The MSE (\ref{mse}) of $\wW_{MMSE}\y$ in (\ref{estimate}) is given under assumption (\ref{asu}) as:
  \begin{multline} \label{mse_bubu_obtain}
    \Jn(\wW_{MMSE})=\sum_{i=1}^m{{\si_i^2(\ga_i^2+\ep)}\over
      {(\si_i^2+\ep)^2}}-
    2\sum_{i=1}^m{{\ga_i\si_i}\over{\si_i^2+\ep}}+m=\\\sum_{i=1}^m x_i^2(\ga_i^2+\ep)-
    2\sum_{i=1}^mx_i\ga_i+m,
  \end{multline}
  where 
  \begin{equation} \label{ixi}
    x_i=\frac{\si_i}{\si_i^2+\ep},\ i=1,\dots,m.
  \end{equation} 
  The MSE (\ref{mse}) of $\wW_{r-MMSE}\y$ in (\ref{estimate}) is given under assumption (\ref{asu}) as:
  \begin{multline} \label{rrmse_bubu_obtain}
    \Jn(\wW_{r-MMSE})=\sum_{i=1}^r{{\si_i^2(\ga_i^2+\ep)}\over
      {(\si_i^2+\ep)^2}}-
    2\sum_{i=1}^r{{\ga_i\si_i}\over{\si_i^2+\ep}}+m=\\\sum_{i=1}^r x_i^2(\ga_i^2+\ep)-
    2\sum_{i=1}^rx_i\ga_i+m.
  \end{multline}
  The MSE (\ref{mse}) of $\wW_{r-SVD}\y$ in (\ref{estimate}) is given under assumption (\ref{asu}) as:
  \begin{equation} \label{mvp_bubu_obtain}
    \Jn(\wW_{r-SVD})=\sum_{i=1}^r\frac{\ga_i^2+\ep}{\si_i^2}-2\sum_{i=1}^r\frac{\ga_i}{\si_i}+m=m-r+\ep\sum_{i=1}^r{{1}\over{\si_i^2}}+\sum_{i=1}^r\biggl({{\ga_i}\over{\si_i}}-1\biggr)^2.
  \end{equation}
\end{proposition}
\pd{See \ref{pd_s_expr}.}

We define now the following functions:
\begin{equation} \label{A_i}
  A_i(\ga_i,\si_i,\ep)=\Biggl[x_i^2(\ga_i^2+\ep)-2x_i\ga_i\Biggr]-
  \left[\ep{{1}\over{\si_i^2}}+
    \left({{\ga_i}\over{\si_i}}-1\right)^2\right]+1,\  i=1,\dots,m,
\end{equation}
and
\begin{equation} \label{B_i}
  B_i(\ga_i,\si_i,\ep)=x_i^2(\ga_i^2+\ep)-
  2x_i\ga_i, \  i=1,\dots,m,
\end{equation}
where $x_i$ is defined as in (\ref{ixi}). Then, it is simple to verify that:
\begin{equation} \label{M}
  \Jn(\wW_{MMSE})=\sum_{i=1}^mB_i(\ga_i,\si_i,\ep)+m,
\end{equation}
\begin{equation} \label{R}
  \Jn(\wW_{r-MMSE})=\sum_{i=1}^rB_i(\ga_i,\si_i,\ep)+m,
\end{equation}
\begin{equation} \label{R-T}
  \Jn(\wW_{r-MMSE})-\Jn(\wW_{r-SVD})=\sum_{i=1}^rA_i(\ga_i,\si_i,\ep),
\end{equation}
\begin{equation} \label{M-R}
  \Jn(\wW_{MMSE})-\Jn(\wW_{r-MMSE})=\sum_{i=r+1}^mB_i(\ga_i,\si_i,\ep),
\end{equation}
and
\begin{equation} \label{M-T}
  \Jn(\wW_{MMSE})-\Jn(\wW_{r-SVD})=\sum_{i=1}^rA_i(\ga_i,\si_i,\ep)+\sum_{i=r+1}^mB_i(\ga_i,\si_i,\ep).
\end{equation}

\begin{proposition} \label{elilama}
  Let $i\in\{1,\dots,m\}.$ Then, for fixed $\si_i$ and $\ep$, $A_i$ in (\ref{A_i}) is a concave function of $\ga_i$ and $B_i$ in (\ref{B_i}) is a convex function of $\ga_i.$ The maximum of $A_i$ is achieved at
  \begin{equation} \label{ingen}
    \ga_i^{A_i^{\max}}={{\si_i(\si_i^2+\ep)}\over{2\si_i^2+\ep}},\ i=1,\dots,m,
  \end{equation}
  while the minimum of $B_i$ is achieved at 
  \begin{equation} \label{thorn}
    \ga_i^{B_i^{\min}}=x_i^{-1}=\si_i+{{\ep}\over{\si_i}},\ i=1,\dots,m.
  \end{equation}
\end{proposition}
\pd See \ref{pd_elilama}. 

\begin{remark}
  If no perturbation is present, i.e., if $\De H=0$, we obviously have $\si_i=\ga_i$ for $i=1,\dots,m.$ Surprisingly, for fixed $\si_i$ and~$\ep$, $\Jn(\wW_{MMSE})$ and $\Jn(\wW_{r-MMSE})$ do not achieve their global minima for $\si_i=\ga_i$ for $i=1,\dots,m.$ Instead, their global minima as functions of $\ga_i$ are defined by the triplets $(x_i^{-1},\si_i,\ep)$ for $i=1,\dots,m$, where $x_i$ is given in (\ref{ixi}).
\end{remark}

Having established concavity of $A_i$ and convexity of $B_i$, we find now intervals in which $A_i$ for $i=1,\dots,r$ and $B_i$ for $i=r+1,\dots,m$ are positive as functions of $\ga_i.$ 

\begin{proposition} \label{wolfram}
  Let $r$ be such that $1\leq r<m.$ Then, for
  fixed $\si_i$ and $\ep$ the following hold:
  \begin{enumerate}
  \item for $i=1,\dots,r$, if $0<\ep<0.325\si_i^2$, then $\be_i^A={\sqrt{{\si_i^6-2\ep\si_i^4-3\ep^2\si_i^2-\ep^3}}\over{\ep+2\si_i^2}}>0$ and
    \begin{equation} \label{whew1eq}
      A_i(\ga_i,\si_i,\ep)>0 \iffm \ga_i\in\left(\ga_i^{A_i^{\max}}-\be_i^A,\ \ga_i^{A_i^{\max}}+\be_i^A\right),
    \end{equation}
    where $\ga_i^{A_i^{\max}}={{\si_i(\si_i^2+\ep)}\over{2\si_i^2+\ep}}$,  \label{whew1}
  \item for $i=r+1,\dots,m$,
    \begin{equation} \label{whew2eq}
      B_i(\ga_i,\si_i,\ep)>0 \iffm \ga_i\in\left(0,\ x_i^{-1}-\be_i^B\right)\cup\left(x_i^{-1}+\be_i^B,\ \infty\right),
    \end{equation}
    where $x_i^{-1}=\si_i+{{\ep}\over{\si_i}}$ and $\be_i^B=\sqrt{\si_i^2+{{\ep^2}\over{\si_i^2}}+\ep}.$ \label{whew2}
  \end{enumerate}
\end{proposition}
\pd See \ref{pd_wolfram}.

A few remarks on Proposition \ref{wolfram} are in place here.

\begin{remark} \label{YS3}
  In view of condition $0<\ep<0.325\si_i^2\iffm\si_i>1.75\sqrt{\ep}$ in \ref{whew1}. in Proposition \ref{wolfram}, we set $\kappa=1.75$ in condition \ref{c2}. in Definition \ref{highill} for the remaining part of this section. Moreover, for $\si_i>1.75\sqrt{\ep}$ one has that $\ga_i^{A_i^{\max}}+\be_i^A<\si_i$ in (\ref{whew1eq}), which can be perhaps most easily verified using Mathematica due to rather complex expression for $\ga_i^{A_i^{\max}}+\be_i^A$ in (\ref{whew1eq}).
\end{remark}

\begin{remark} \label{YS2}
  Since $x_i^{-1}+\be_i^B>\sqrt{\ep}$, only the first interval in (\ref{whew2eq}) is of interest, as $\ga_i\in\left(x_i^{-1}+\be_i^B,\ \infty\right)$ for $i=r+1,\dots,m$ would violate condition \ref{cend}. in Definition \ref{highill}. Moreover, from convexity of $B_i$ as a function of $\ga_i$ only, from (\ref{thorn}) and (\ref{whew2eq}) we conclude that any maximum of $B_i$ as a function of $\ga_i$ and $\si_i$ (i.e., for fixed $\ep$) for $\ga_i\in\left(0,\ x_i^{-1}-\be_i^B\right)$  may be achieved only for $\ga_i\to 0.$
\end{remark}

The following proposition introduces parametrization of $\ga_i$ and $\si_i$ which enables us to eliminate $\ep$ from (\ref{A_i}) and~(\ref{B_i}), and thus gives the same interpretation of values of $A_i$ and $B_i$ for any given noise power~$\ep.$ 
\begin{proposition} \label{wolfram2}
  Let $r$ be such that $1\leq r<m$, and consider $A_i$ in (\ref{A_i}) for $i\in\{1,\dots,r\}$ and $B_i$ in (\ref{B_i}) for $i\in\{r+1,\dots,m\}.$ For fixed $\ep$, let us express values of $\ga_i$ and $\si_i$ as $\ga_i=a_{\ga_i}\sqrt{\ep}$ and $\si_i=a_{\si_i}\sqrt{\ep}$ for some $a_{\ga_i}>0$ and $a_{\si_i}>0$, respectively. Then:
  \begin{equation} \label{ja}
    A_i(a_{\ga_i}\sqrt{\ep},a_{\si_i}\sqrt{\ep},\ep)={{2a_{\ga_i}a_{\si_i}(1+a_{\si_i}^2)-(1+2a_{\si_i}^2)(1+2a_{\ga_i}^2)}\over{(a_{\si_i}+a_{\si_i}^3)^2}},\ i=1,\dots,r,
  \end{equation}
  and
  \begin{equation} \label{pan}
    B_i(a_{\ga_i}\sqrt{\ep},a_{\si_i}\sqrt{\ep},\ep)={{a_{\si_i}[a_{\si_i}+a_{\ga_i}(a_{\si_i}a_{\ga_i}-2a_{\si_i}^2-2)]}\over{(1+a_{\si_i}^2)^2}},\ i=r+1,\dots,m.
  \end{equation}
  The unique maximum of $A_i$ for fixed $\ep$ is of the following form:
  \begin{equation} \label{kamieni}
    A_i(1.394\sqrt{\ep},2.611\sqrt{\ep},\ep)\approx 0.033,\ i=1,\dots,r.
  \end{equation}
  Let us now restrict the domain of $\ga_i$ to the interval $\left(0,\ x_i^{-1}-\be_i^B\right)$, where $x_i^{-1}=\si_i+{{\ep}\over{\si_i}}$ and $\be_i^B=\sqrt{\si_i^2+{{\ep^2}\over{\si_i^2}}+\ep}$, \emph{cf.} Remark \ref{YS2}. Then, for $\si_i=\sqrt{\ep}$ and fixed $\ep$ one has:
  \begin{equation} \label{kupa}
    \lim_{\ga_i\rightarrow 0}B_i(\ga_i,\sqrt{\ep},\ep)=0.25,\ i=r+1,\dots,m,
  \end{equation}
  which is the supremum of $B_i$ for $\ga_i\in\left(0,\ x_i^{-1}-\be_i^B\right)$, $\si_i>0$, and fixed $\ep>0.$
\end{proposition} 
\pd See \ref{pd_wolfram2}. 

\begin{remark}
  In view of Remark \ref{YS3} below Proposition \ref{wolfram}, the value of $\ga_i=1.394\sqrt{\ep}$ in (\ref{kamieni}) does not satisfy the condition $\ga_i>\kappa\sqrt{\ep}$ in Definition \ref{highill} for $\kappa=1.75.$
\end{remark}

\subsection{Interpretation of Propositions \ref{wolfram} and \ref{wolfram2}} \label{interpreter}
We illustrate in Figs.\ref{fig:1}-\ref{fig:4} the results of Propositions \ref{wolfram} and \ref{wolfram2}, where we plot $A_i(a_{\ga_i}\sqrt{\ep},a_{\si_i}\sqrt{\ep},\ep)$ expressed as in (\ref{ja}) and $B_i(a_{\ga_i}\sqrt{\ep},a_{\si_i}\sqrt{\ep},\ep)$ expressed as in (\ref{pan}) as functions of $a_{\ga_i}>0$ and $a_{\si_i}>0$ with noise power $\ep$ as a fixed parameter. We also describe below certain sufficient conditions derived from the results of Propositions~\ref{wolfram}~and~\ref{wolfram2}. 

\begin{enumerate}
\item $\Jn(\wW_{r-MMSE})-\Jn(\wW_{r-SVD})=\sum_{i=1}^rA_i(a_{\ga_i}\sqrt{\ep},a_{\si_i}\sqrt{\ep},\ep)$ in (\ref{R-T}): if the leading $r$ singular values of $H$, $\ga_1,\dots,\ga_r$, are within the interval $\ga_i\in\left(\ga_i^{A_i^{\max}}-\be_i^A,\ \ga_i^{A_i^{\max}}+\be_i^A\right)$ in (\ref{whew1eq}), then $r$-SVD estimator will achieve lower MSE than $r$-MMSE estimator. We also note that from (\ref{ja}) it is seen that $\lim_{\si_i\to\infty}A_i(\ga_i,\si_i,\ep)=0.$ \label{iR-T}
\item $\Jn(\wW_{MMSE})-\Jn(\wW_{r-MMSE})=\sum_{i=r+1}^mB_i(a_{\ga_i}\sqrt{\ep},a_{\si_i}\sqrt{\ep},\ep)$ in (\ref{M-R}): if the trailing $m-r$ singular values of $H$ are such that $\ga_i\in\left(0,\ x_i^{-1}-\be_i^B\right)$ in (\ref{whew2eq}), then $r$-MMSE estimator will achieve lower MSE than MMSE estimator. We note that in view of Remark \ref{YS2} we do not consider $\ga_i\in\left(x_i^{-1}+\be_i^B,\ \infty\right)$ for $i=r+1,\dots,m.$

  \pagebreak
  \begin{figure}[!ht] 
    \centering
    \includegraphics[scale=.45]{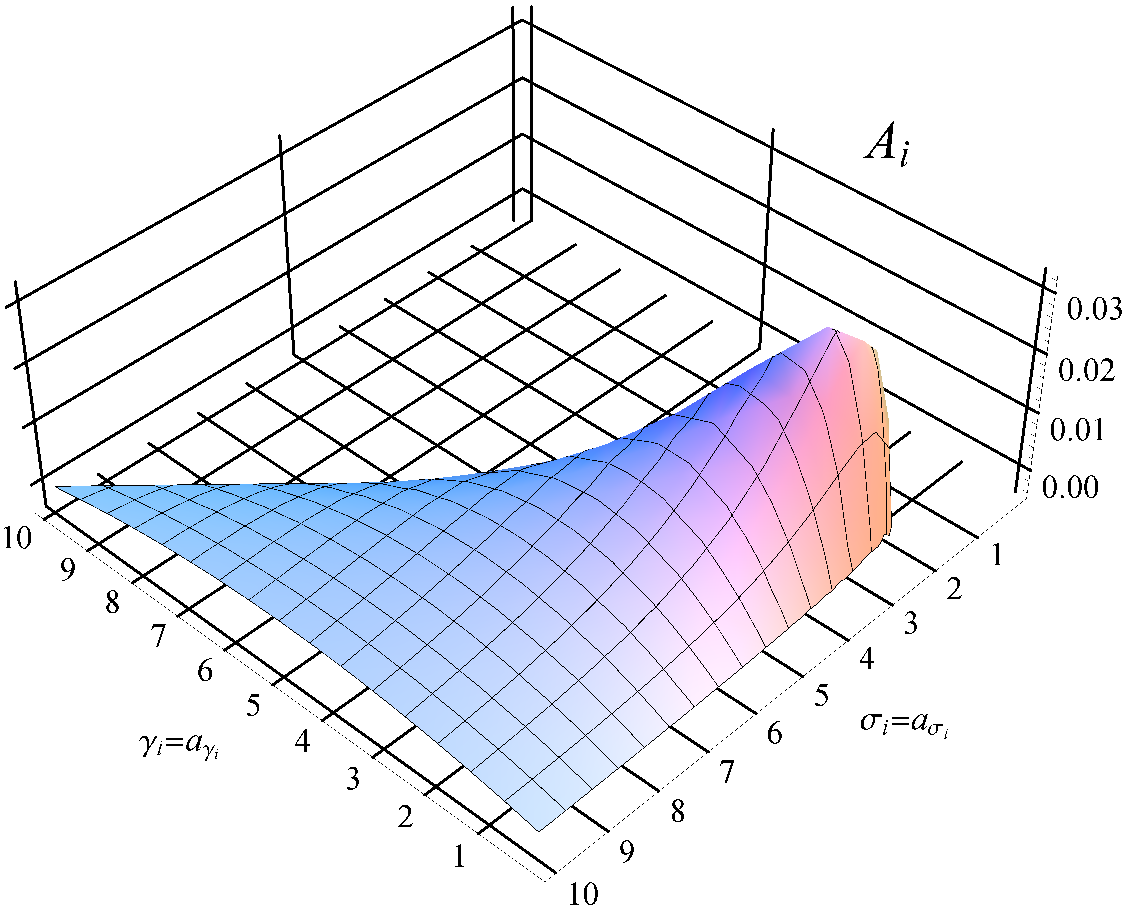} 
    \caption{$A_i$ evaluated over $(a_{\ga_i},a_{\si_i})$ domain yielding positive values, for any given noise level $\ep.$}
    \label{fig:1}
  \end{figure}

  \begin{figure}[!hb] 
    \centering
    \includegraphics[scale=.4]{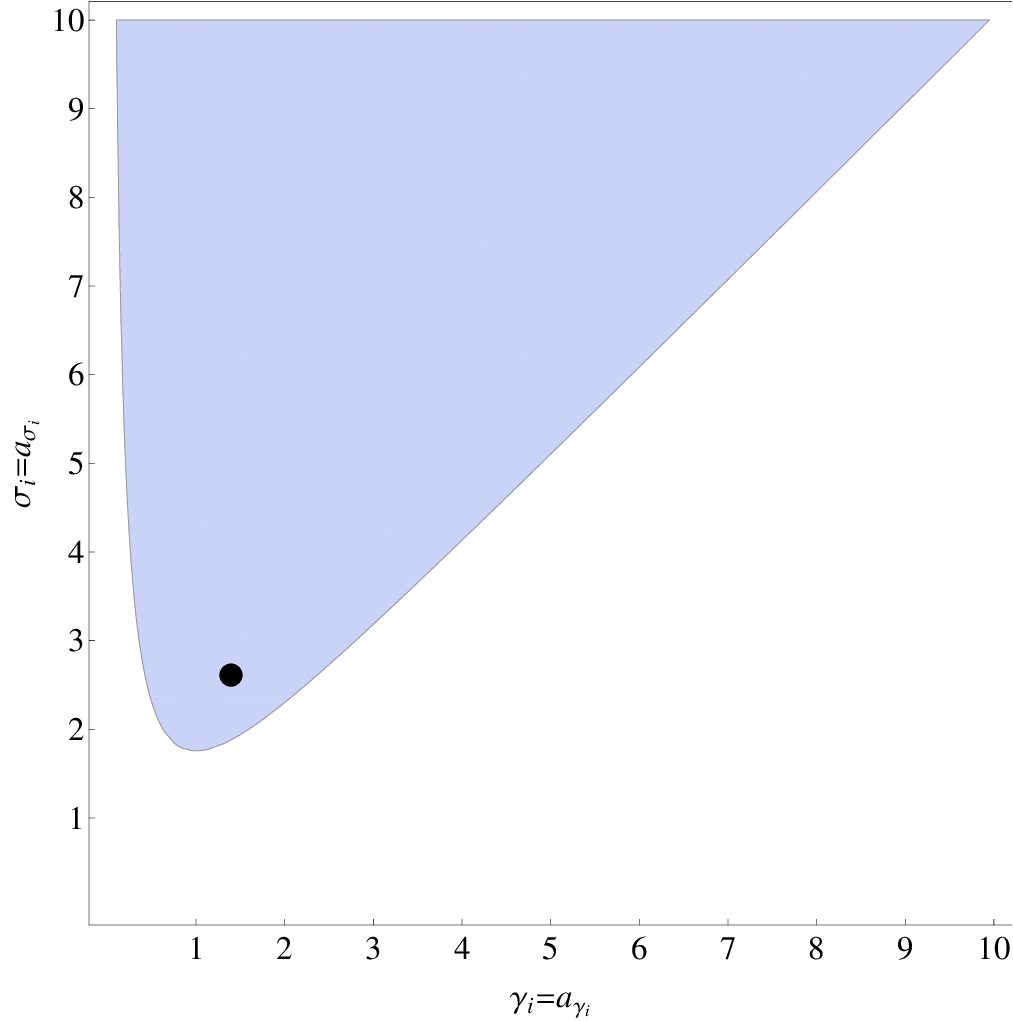} 
    \caption{Domain of $(a_{\ga_i},a_{\si_i})$ yielding positive values of $A_i$, for any given noise level $\ep.$ The unique argument of maximum in (\ref{kamieni}) is marked by a black~dot.}
    \label{fig:2}
  \end{figure}
  \pagebreak

  \begin{figure}[!ht] 
    \centering
    \includegraphics[scale=.45]{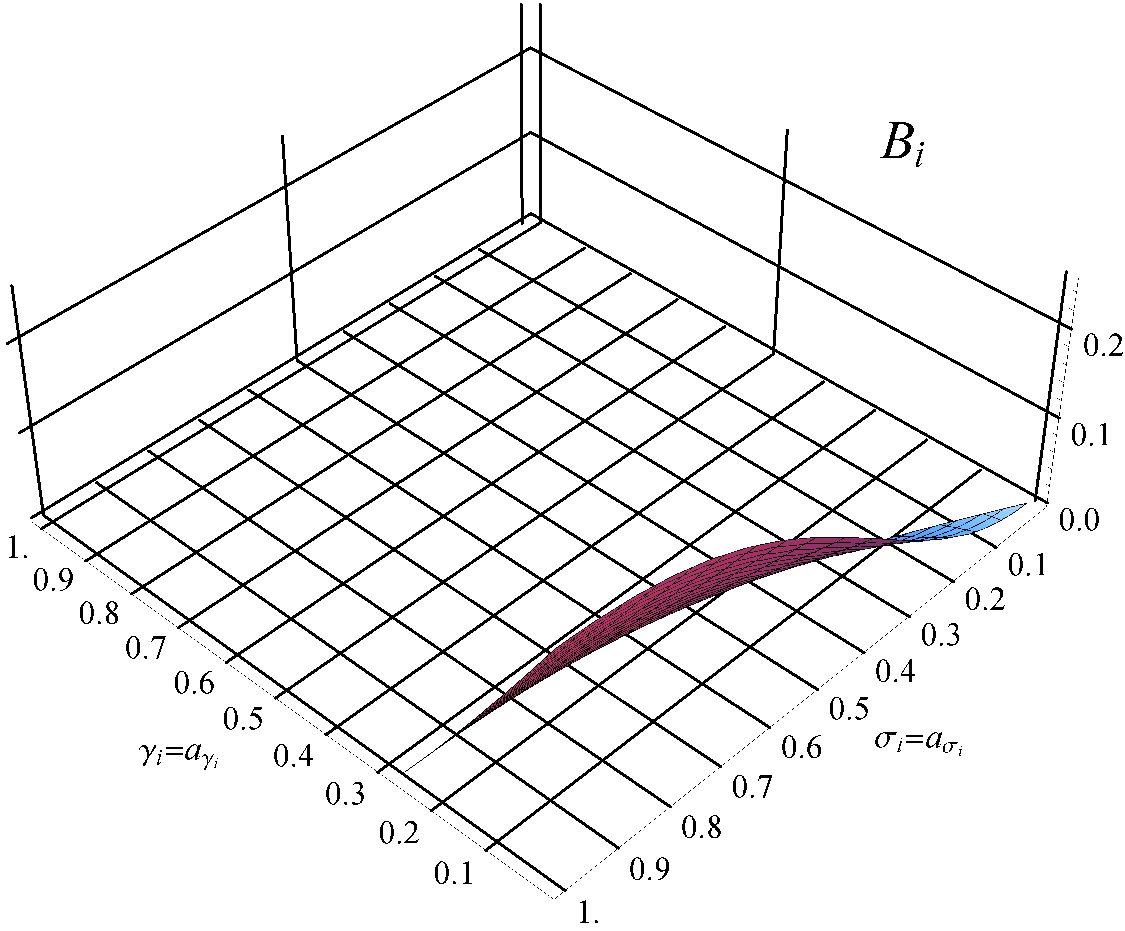} 
    \caption{$B_i$ evaluated over $(a_{\ga_i},a_{\si_i})$ domain yielding positive values, for any given noise level $\ep.$}
    \label{fig:3}
  \end{figure}

  \begin{figure}[!hb] 
    \centering
    \includegraphics[scale=.4]{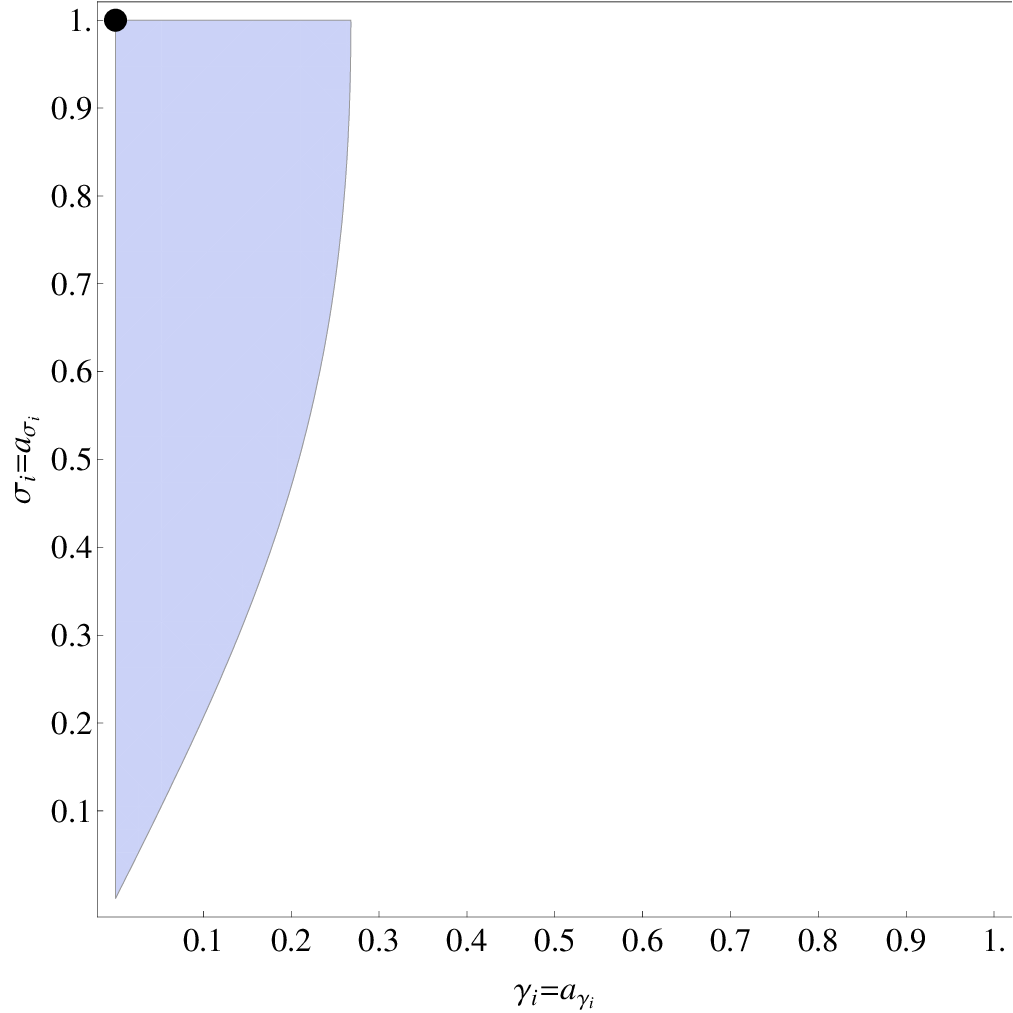} 
    \caption{Domain of $(a_{\ga_i},a_{\si_i})$ yielding positive values of $B_i$, for any given noise level $\ep.$ The unique argument of supremum in (\ref{kupa}) is marked by a black~dot.}
    \label{fig:4}
  \end{figure}
  \pagebreak
  
  We also note that since $\si_i\leq\sqrt{\ep}$ for $i=r+1,\dots,m$, then $\ga_i<(2-\sqrt{3})\sqrt{\ep}\approx 0.27\sqrt{\ep}$ in (\ref{whew2eq}) in view of $\ga_i\in\left(0,\ x_i^{-1}-\be_i^B\right)$ for $i=r+1,\dots,m.$ Hence, $H$ has vanishingly small singular values $\{\ga_{r+1},\dots,\ga_m\}.$ In such a case, the set $\{\si_{r+1},\dots,\si_m\}$ is likely to contain larger singular values than $\{\ga_{r+1},\dots,\ga_m\}$, as small singular values tend to get larger under perturbation, see \cite[pp.263-269]{Stewart1990} and \cite[p. 216]{Stewart1979} for details. Thus, for $\ga_i\in\left(0,\ x_i^{-1}-\be_i^B\right)$ for $i=r+1,\dots,m$ we are likely to obtain distribution of $\{\si_{r+1},\dots,\si_m\}$ yielding positive values of $B_i$, as depicted in Fig.\ref{fig:4}. We will also observe this fact in the numerical example in Section \ref{ne}. 
  \label{iM-R}
\item $\Jn(\wW_{MMSE})-\Jn(\wW_{r-SVD})=\\\sum_{i=1}^rA_i(a_{\ga_i}\sqrt{\ep},a_{\si_i}\sqrt{\ep},\ep)+\sum_{i=r+1}^mB_i(a_{\ga_i}\sqrt{\ep},a_{\si_i}\sqrt{\ep},\ep)$ in (\ref{M-T}): this expression is a sum of expressions in \ref{iR-T}) and \ref{iM-R}), thus discussion of both of the above points applies. \label{iM-T}
\end{enumerate}

\section{Extension to Generic Perturbations} \label{ext}
We drop now assumption (\ref{asu}), and consider the case of an arbitrary perturbation $\De H$ such that Definition \ref{highill} remains valid.

\subsection{Objectives} \label{obj}
We show that for $r$ such that $1\leq r<m$ and a given distribution of singular values of $H$ and $H+\De H$ satisfying Definition \ref{highill}, the conclusions of Section \ref{simple} remain approximately valid. More precisely, we show that:
\begin{enumerate}
\item the values of $\Jn(\wW_{r-MMSE})$ and $\Jn(\wW_{r-SVD})$ are, respectively, approximately equal to (\ref{rrmse_bubu_obtain}) and (\ref{mvp_bubu_obtain}) obtained under assumption (\ref{asu}),
  \label{Viv1}
\item the value of $\Jn(\wW_{MMSE})$ is likely to be larger than $\Jn(\wW_{r-MMSE}).$ 
  \label{Viv2}
\end{enumerate}

\subsection{Results} \label{res}
We note first that $\Ry$ can be expressed in terms of the singular value decomposition of $H$ in (\ref{svdtruH}) as 

\begin{equation} \label{no1}
  \Ry=HH^t+\ep I_n=M(\Ga\Ga^t+\ep I_n)M^t,
\end{equation}
and that $\Ryx=H=M\Ga N^t.$ We also define [\emph{cf.} (\ref{svdtruH}) and (\ref{svdH})]:
\begin{equation} \label{X}
  X=\Si^t(\Si\Si^t+\ep I_n)^{-1}\in\mn,
\end{equation}
\begin{equation} \label{Y}
  Y=\Ga\Ga^t+\ep I_n\in\nn,
\end{equation}
\begin{equation} \label{K}
  K=U^tM\in\nn,
\end{equation}
\begin{equation} \label{L}
  L=N^tV\in\mm,
\end{equation}
with $k_{i,j}$ for $i=1,\dots,n$ and $j=1,\dots,n$ denoting the $i,j$-th entry of $K$, and similarly $l_{i,j}$ for $i=1,\dots,m$ and $j=1,\dots,m$ denoting the $i,j$-th entry of $L.$ We note that both $K$ and $L$ are orthogonal matrices.

\begin{proposition} \label{mse_exprs}
  Consider MMSE (\ref{mmseish}), $r$-MMSE (\ref{rrmmseish}), and $r$-SVD (\ref{mvpform_s}) estimators designed based on model (\ref{bubu}), and let us denote
  \begin{equation} \label{phipsi}
    \phi_i=\sqrt{\sum_{j=1}^mk_{i,j}^2\ga_j^2},\ 
    \psi_i=\sum_{j=1}^m\ga_jk_{i,j}l_{j,i},\ i=1,\dots,m.
  \end{equation}
  Then, the MSE (\ref{mse}) of $\wW_{MMSE}\y$ in (\ref{estimate}) is given as:
  \begin{equation} \label{sami}
    \Jn(\wW_{MMSE})=
    \sum_{i=1}^mx_i^2(\phi_i^2+\ep)-2\sum_{i=1}^mx_i\psi_i+m,
  \end{equation}
  where $x_i=\frac{\si_i}{\si_i^2+\ep}$ is defined in (\ref{ixi}). Similarly, for a given rank constraint $r$ such that $1\leq r<m$, the MSE (\ref{mse}) of $\wW_{r-MMSE}\y$ in (\ref{estimate}) is given as:
  \begin{equation} \label{samija}
    \Jn(\wW_{r-MMSE})=
    \sum_{i=1}^rx_i^2(\phi_i^2+\ep)-2\sum_{i=1}^rx_i\psi_i+m,
  \end{equation}
  while the MSE (\ref{mse}) of $\wW_{r-SVD}\y$ in (\ref{estimate}) is given as:
  \begin{equation} \label{kedira}
    \Jn(\wW_{r-SVD})=\sum_{i=1}^r\frac{\phi_i^2+\ep}{\si_i^2}-2\sum_{i=1}^r\frac{\psi_i}{\si_i}+m.
  \end{equation}
\end{proposition}
\pd See \ref{pd_mse_exprs}.

\begin{corollary} \label{C_main}
  Let us compare now the MSE expressions of $\wW_{r-MMSE}\y$ and $\wW_{r-SVD}\y$ obtained under assumption (\ref{asu}) in Proposition \ref{s_expr} in Section~\ref{simple} with those obtained in Proposition \ref{mse_exprs} above for the generic case. For the sake of this comparison, we denote the former expressions as $_{(\ref{simple})}\Jn(\wW_{r-MMSE})$ and $_{(\ref{simple})}\Jn(\wW_{r-SVD})$, and the latter as
  $_{(\ref{ext})}\Jn(\wW_{r-MMSE})$ and $_{(\ref{ext})}\Jn(\wW_{r-SVD})$, respectively. Then:
  \begin{equation} \label{el}
    \left\lvert_{(\ref{ext})}\Jn(\wW_{r-MMSE})-_{(\ref{simple})}\Jn(\wW_{r-MMSE})\right\rvert=\left\lvert
    \sum_{i=1}^r\left[x_i^2\left(\phi_i^2-\ga_i^2\right)+
      2x_i\left(\ga_i-\psi_i\right)\right]\right\rvert,
  \end{equation}
  and
  \begin{equation} \label{ter}
    \left\lvert_{(\ref{ext})}\Jn(\wW_{r-SVD})-_{(\ref{simple})}\Jn(\wW_{r-SVD})\right\rvert=\left\lvert\sum_{i=1}^r\left[\si_i^{-2}\left(\phi_i^2-\ga_i^2\right)+2\si_i^{-1}\left(\ga_i-\psi_i\right)\right]\right\rvert.
  \end{equation}
  We also have:
  \begin{equation} \label{al}
    _{(\ref{ext})}\Jn(\wW_{MMSE})-_{(\ref{ext})}\Jn(\wW_{r-MMSE})=
    \sum_{i=r+1}^mx_i^2(\phi_i^2+\ep)-2\sum_{i=r+1}^mx_i\psi_i,
  \end{equation}
  where $_{(\ref{ext})}\Jn(\wW_{MMSE})$ is the MSE expression given in (\ref{sami}).
\end{corollary}

The following lemma establishes bounds on the first $r$ diagonal entries of $K$ in (\ref{K}) and $L$ in (\ref{L}). In particular, the term $\rho_i$ in (\ref{rho}) will play a pivotal role in subsequent derivations.
\begin{lemma} \label{fina}
  Let us consider SVD of $H$ of the form (\ref{svdtruH}) and SVD of $H+\De H$ of the form (\ref{svdH}), and let $r$ such that $1\leq r<m$ be selected according to Definition \ref{highill}. For $i=1,\dots,r$  we denote:
  \begin{equation} \label{rho}
    \rho_i=1-{{\p\De H\p_2^2}\over{\de_i^2}},
  \end{equation}
  where 
  \begin{equation} \label{mex}
    \de_i=\min_{{\ga}\in\bs{\ga_{i, NEXT}}}|\si_i-\ga|,
  \end{equation}
  with
  \begin{equation} \label{27}
    \bs{\ga_{1, NEXT}}=\{\ga_2\},
  \end{equation}
  \begin{equation} \label{08}
    \bs{\ga_{i,NEXT}}=\{\ga_{i-1},\ga_{i+1}\},\ i\in\{2,\dots,r\}.
  \end{equation}
  Then $k_{i,i}$ and $l_{i,i}$ in (\ref{phipsi})-(\ref{kedira}) satisfy: 
  \begin{equation} \label{faaaa}
    2\geq\lvert k_{i,i}\rvert+\lvert l_{i,i}\rvert\geq k_{i,i}^2+l_{i,i}^2\geq 2\rho_i.
  \end{equation}
  Moreover, if $\rho_i>1/2$, then
  \begin{equation} \label{fuuuu}
    1\geq k_{i,i}l_{i,i}>0,
  \end{equation}
  i.e., $k_{i,i}$ and $l_{i,i}$ have matching signs for $i=1,\dots,r.$
\end{lemma}
\pd See \ref{pd_fina}.

Lemma \ref{fina} allows one to prove the following proposition.
\begin{proposition} \label{funda}
  With notation as in Lemma \ref{fina}, consider $i\in\{1,\dots,r\}$ such that $r<m$, and assume that $\rho_i>1/2.$ The following inequalities hold:
  \begin{equation} \label{conc0}
    k_{i,i}^2\geq 2\rho_i-1,\ l_{i,i}^2\geq 2\rho_i-1,\  k_{i,i}l_{i,i}\geq 2\rho_i-1,
  \end{equation}
  \begin{multline} \label{conc1}
    \sum_{j=1,j\neq i}^nk_{i,j}^2\leq 2(1-\rho_i),\ \sum_{j=1,j\neq i}^nk_{j,i}^2\leq 2(1-\rho_i),\\ \zt{ for } i\neq j\zt{ and }j\leq n\ \ k_{i,j}^2\leq 2(1-\rho_i),\ k_{j,i}^2\leq 2(1-\rho_i),
  \end{multline}
  \begin{multline} \label{conc1.5}
    \sum_{j=1,j\neq i}^ml_{i,j}^2\leq 2(1-\rho_i),\ \sum_{j=1,j\neq i}^ml_{j,i}^2\leq 2(1-\rho_i),\\ \zt{ for } i\neq j\zt{ and }j\leq m\ \ l_{i,j}^2\leq 2(1-\rho_i), \ l_{j,i}^2\leq 2(1-\rho_i),
  \end{multline}
  \begin{equation} \label{conc2}
    \left|\sum_{j=1,j\neq i}^mk_{i,j}l_{j,i}\right|\leq 2(1-\rho_i), \zt{ for } i\neq j\zt{ and }j\leq m\ \ |k_{i,j}l_{j,i}|\leq 2(1-\rho_i).
  \end{equation}
\end{proposition}
\pd See \ref{pd_funda}.

\begin{corollary} \label{yea}
  With notation as in Lemma \ref{fina} and for $\rho_i$ in (\ref{rho}) close to 1 for $i=1,\dots,r$, the orthogonal matrices $K$ in (\ref{K}) and $L$~in~(\ref{L}) can be approximated as:
  \begin{equation} \label{KandL}
    K\approx\left( 
    \begin{array}{cc}
      I^{\pm 1}_r & 0_{r\times (n-r)}\\
      0_{(n-r)\times r} & Z^K_{n-r}\\ 
    \end{array} \right)
    \quad\zt{and}\quad
    L\approx\left( 
    \begin{array}{cc}
      I^{\pm 1}_r & 0_{r\times (m-r)}\\
      0_{(m-r)\times r} & Z^L_{m-r}\\ 
    \end{array} \right), 
  \end{equation}
  where $I^{\pm 1}_r$ is a diagonal matrix of size $r$ with $\pm 1$ on the diagonal, $0_{x\times y}$ is a matrix of zeros of size $x\times y$, and $Z^K_{n-r}$ and $Z^L_{m-r}$ are certain orthogonal matrices of respective sizes. 
\end{corollary}

The following theorem allows us to justify claims \ref{Viv1}. and \ref{Viv2}. in Section~\ref{obj}.
\begin{theorem} \label{greatinequalities}
  Let SVD of $H$ be of the form (\ref{svdtruH}) and SVD of $H+\De H$ of the form (\ref{svdH}), and let $r$ be such that $1\leq r<m$ be selected according to Definition~\ref{highill}. Let $\rho_i$ for $i=1,\dots,r$ be defined as in Lemma \ref{fina} through (\ref{rho})-(\ref{08}), and assume that $\rho_i>1/2$ for $i=1,\dots,r.$ Then, with notation of Corollary \ref{C_main} the following inequalities hold:
  \begin{equation} \label{elGI}
    \left\lvert_{(\ref{ext})}\Jn(\wW_{r-MMSE})-_{(\ref{simple})}\Jn(\wW_{r-MMSE})\right\rvert\leq
    2\sum_{i=1}^r\sum_{j=1}^m(1-\rho_i)(x_i\ga_j+1)^2,
  \end{equation}
  where $x_i=\frac{\si_i}{\si_i^2+\ep}$ has been introduced in (\ref{ixi}). Furthermore,
  \begin{equation} \label{terGI}
    \left\lvert_{(\ref{ext})}\Jn(\wW_{r-SVD})-_{(\ref{simple})}\Jn(\wW_{r-SVD})\right\rvert\leq
    2\sum_{i=1}^r\sum_{j=1}^m(1-\rho_i)(\si_i^{-1}\ga_j+1)^2.
  \end{equation}
  Moreover, $_{(\ref{ext})}\Jn(\wW_{MMSE})\geq_{(\ref{ext})}\Jn(\wW_{r-MMSE})$ if
  \begin{equation} \label{alGI}
    \sum_{i=r+1}^m\sum_{j=1}^m(x_i\ga_jk_{i,j}-l_{j,i})^2+\ep\sum_{i=r+1}^mx_i^2\geq 
    m-r,
  \end{equation}
  where $k_{i,j}$ and $l_{j,i}$ are the entries of matrices $K$ in (\ref{K}) and $L$ in (\ref{L}), respectively.
\end{theorem}
\pd See \ref{pd_greatinequalities}.

\begin{remark} \label{R0}
  In Remarks \ref{R1} and \ref{R2} below we analyze the conditions under which the right-hand sides of inequalities (\ref{elGI})-(\ref{terGI}) are small, and the left-hand side of inequality (\ref{alGI}) may be larger than $m-r.$ Namely, as will be seen in Remarks \ref{R1} and \ref{R2}, the inequalities (\ref{elGI})-(\ref{alGI}) are the most informative if $\rho_i$ in (\ref{rho}) is close enough to 1 for $i=1,\dots,r$ to ensure the existence of a small constant $c>0$ such that
  \begin{equation} \label{ost}
    {{\p\De H\p_2^2}\over{\de_i^2}}\leq c<1/2,\ i=1,\dots,r.
  \end{equation}
  From (\ref{ost}) and (\ref{mex})-(\ref{08}) in Lemma \ref{fina} it is seen that $c$ can be made arbitrarily small by either decreasing the size of perturbation $\De H$ or increasing the spread of the leading $r$ singular values of $H$ and $H+\De H$ and the value of $\kappa$ in condition~\ref{c2}. in Definition \ref{highill}, which in view of (\ref{mex})-(\ref{08}) will increase the value of $\de_i$ for $i=1,\dots,r.$ In Subsection \ref{cs} we discuss when such conditions are met for a random square Gaussian distributed $H.$  Furthermore, based on the analysis of Subsection \ref{cs}, in Section \ref{ne} we present a simple numerical setup where $c=0.0013.$ We note that for such small $c$ [and hence for $\rho_i$ in (\ref{rho}) close to 1] the approximation in Corollary \ref{yea} will be very accurate, as established by inequalities in Proposition \ref{funda}.

  We also assume that $\ga_i$ are vanishingly small for $i=r+1,\dots,m$, and that at least some $\si_i$ are close to $\sqrt{\ep}$ for $i\in\{r+1,\dots,m\}$. We note that these are precisely the conditions which yield the greatest gain in MSE performance of reduced-rank estimators over the MMSE estimator for the perturbed singular values model considered in Section \ref{simple}, as given in Proposition~\ref{wolfram2}, \emph{cf.} also Figs.\ref{fig:3}-\ref{fig:4}.
\end{remark}

\begin{remark}[on inequalities (\ref{elGI}) and (\ref{terGI})] \label{R1}
  We have: 
  \begin{enumerate} 
  \item From the  definition of $\rho_i$ in (\ref{rho}) in Lemma \ref{fina} for $i=1,\dots,r$, it is clear that the smaller $c$ can be selected in (\ref{ost}), the smaller the upper bounds on the right-hand sides of (\ref{elGI}) and (\ref{terGI}) become. 
  \item For $i=1,\dots,r$ and $j=1,\dots,r$, the terms $x_i=\frac{\si_i}{\si_i^2+\ep}<\si_i^{-1}$ and $\si_i^{-1}$ offset the values of $\ga_j$ for $i,j=1,\dots,r$ in $x_i\ga_j\leq x_r\ga_1$ and $\si_i^{-1}\ga_j\leq \si_r^{-1}\ga_1$ on the right-hand sides of (\ref{elGI}) and (\ref{terGI}), respectively. 
  \item For $i=1,\dots,r$ and $j=r+1,\dots,m$, the terms $x_i\ga_j$ in (\ref{elGI}) and $\si_i^{-1}\ga_j$ in (\ref{terGI}) contain the vanishingly small $\ga_j$ for $j=r+1,\dots,m$ under the assumption made in Remark \ref{R0}.\footnote{However, we note that this assumption is not necessary to keep $x_i\ga_j$ and $\si_i^{-1}\ga_j$ small, as one has $\ga_j\leq\sqrt{\ep}$ for $j=r+1,\dots,m$ from condition~\ref{cend}. in Definition \ref{highill}.}
  \end{enumerate}
\end{remark}

\begin{remark}[on inequality (\ref{alGI})] \label{R2}
  We have:
  \begin{enumerate}
  \item The term $\sum_{i=r+1}^m\sum_{j=1}^r(x_i\ga_jk_{i,j}-l_{j,i})^2$: from (\ref{conc1})-(\ref{conc1.5}) in Proposition~\ref{funda} it is seen in particular that $k_{i,j}^2\leq 2(1-\rho_j)$ and $l_{j,i}^2\leq 2(1-\rho_j)$ for $i=r+1,\dots,m$ and $j=1,\dots,r$ in (\ref{alGI}), where $k_{i,j}$ and $l_{j,i}$ are the entries of matrices $K$ in (\ref{K}) and $L$ in (\ref{L}), respectively.\footnote{We changed the index from $i$ to $j$ for $\rho_j$ in this remark, as it aligns with the scope of indices $i$ and $j$ considered in inequality (\ref{alGI}).} Thus, $|k_{i,j}|\leq\sqrt{2(1-\rho_j)}$ and $|l_{j,i}|\leq\sqrt{2(1-\rho_j)}$ will be vanishingly small for $i=r+1,\dots,m$ and $j=1,\dots,r.$

    However, we note that if $\si_i$ is close to $\sqrt{\ep}$ for $i\in\{r+1,\dots,m\}$ (\emph{cf.} Figs.\ref{fig:3}-\ref{fig:4}), the value of $x_i=\frac{\si_i}{\si_i^2+\ep}$ [see (\ref{ixi})] approaches its unique maximum $\frac{1}{2\sqrt{\ep}}$ attainable for $\si_i=\sqrt{\ep}.$ Thus, for small values of $\ep$ (satisfying Definition~\ref{highill}), the value of $x_i$ may be large for $i\in\{r+1,\dots,m\}.$ Moreover, for a given $\ep$, the larger the value of $\kappa$ in condition~\ref{c2}. in Definition~\ref{highill}, the larger the values of $\ga_j$ for $j=1,\dots,r.$ Thus, from the above analysis we conclude that the value of $\sum_{i=r+1}^m\sum_{j=1}^r(x_i\ga_jk_{i,j}-l_{j,i})^2$ on the left-hand side of (\ref{alGI}) may be nonvanishing positive due to some large $x_i$ and large $\ga_j$ for $i=r+1,\dots,m$ and $j=1,\dots,r.$ \label{Mi}
  \item The term $\sum_{i=r+1}^m\sum_{j=r+1}^m(x_i\ga_jk_{i,j}-l_{j,i})^2$: we note first that in this case $k_{i,j}$ and $l_{j,i}$ are approximately the entries of orthogonal matrices $Z^K_{n-r}$
    \pagebreak

    and $Z^L_{m-r}$ introduced in Corollary \ref{yea}.\footnote{The quality of this approximation is controlled by inequalities (\ref{conc0})-(\ref{conc1.5}) in Proposition~\ref{funda}.} Then, for $\ga_j$ significantly smaller than $\si_i\leq\sqrt{\ep}$ for $i,j=r+1,\dots,m$, the terms $x_i\ga_jk_{i,j}$ will be vanishingly small, and hence 
    $\sum_{i=r+1}^m\sum_{j=r+1}^m(x_i\ga_jk_{i,j}-l_{j,i})^2$ will be close to $\sum_{i=r+1}^m\sum_{j=r+1}^m(-l_{j,i})^2\approx 
    tr[(Z^L_{m-r})^t(Z^L_{m-r})]=tr[I_{m-r}]$ $=m-r.$ 
  \item Finally, we note that for a given noise power $\ep$, the term $\ep\sum_{i=r+1}^mx_i^2$ achieves its unique maximum for $\si_i=\sqrt{\ep}$ for $i=r+1,\dots,m$, with the maximum value $\ep\sum_{i=r+1}^m\frac{1}{4\ep}=\frac{1}{4}(m-r).$
    \label{sc}
  \end{enumerate}
\end{remark}

We will present a simple numerical setup in Section \ref{ne}, where the inequalities in (\ref{elGI})-(\ref{alGI}) will be used simultaneously.

\subsection{Case Study: Square Gaussian Distributed $H$} \label{cs}
As discussed in Remarks \ref{R1}-\ref{R2}, the results of Theorem \ref{greatinequalities} are most useful if $\rho_i$ in (\ref{rho}) is close to $1$ for $i=1,\dots,r.$ In this section we discuss their application to the case of Gaussian distributed $H$ in model (\ref{general}). For simplicity, we focus here on the square case $m=n.$ 

Namely, assume that the entries of $H$ are independent standard normal random variables.\footnote{By simple rescaling, the above analysis is also applicable if the standard deviation of elements of $H$ is different than one.} Then \cite[Th.2.6]{Rudelson2010}
\begin{equation}
  \sqrt{n}-\sqrt{m}\leq\E[\ga_m],
\end{equation}
where $E[\ga_m]$ is the expected value of $\ga_m.$ We note that the lower bound on $\E[\ga_m]$ becomes zero for square matrices. Indeed, the following result gives in particular the probability that $\ga_m$ is smaller than any fixed positive constant for $m=n$ \cite[Th.1.1]{Tao2010}:
\begin{equation} \label{Pre}
  \Pre(m\ga_m^2\leq t)=1-\exp(-t/2-\sqrt{t})+o(1),\ t\geq 0,
\end{equation}
where $\Pre$ is the probability measure. On the other hand, it can be easily seen that one has $\E[\p H\p_2]=\E[\ga_1]\geq\sqrt{\max\{m,n\}}$, see, e.g., \cite{Hansen1988}. Thus, from this fact and (\ref{Pre}) we may assume that in the square-case $H$ may be arbitrarily ill-conditioned with non-zero probability. This hypothesis is confirmed by the following result \cite[Th.1.1]{Azais2004} for $m=n\geq 3$:
\begin{equation} \label{miles}
  \frac{d}{t}<\Pre(cond(H)>mt)<\frac{D}{t},\ t>0,
\end{equation} 
where $cond(H)=\ga_1/\ga_m$ is the condition number of $H$, $d=0.13$ and $D=5.60$ satisfy (\ref{miles}) for any $m=3,4,\dots$ and $t>0.$ Moreover, simulations in \cite{Azais2004} showed that it must be $d<0.87$ and $D>2.18.$ 

\begin{figure}[!ht] 
    \centering
    \includegraphics[scale=.3]{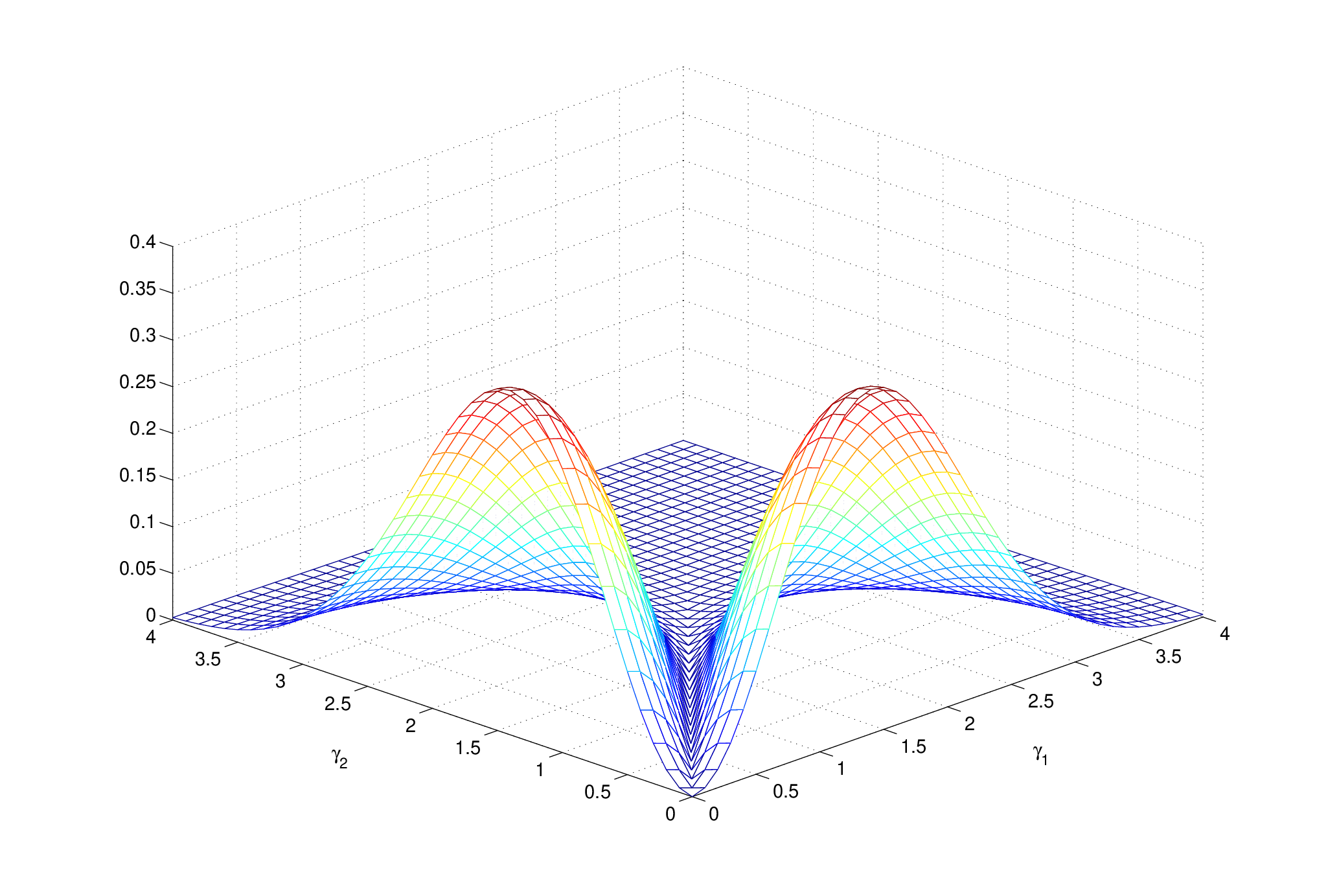} 
    \caption{Probability density of the distribution of singular values of a $2\times 2$ matrix H of independent standard normal random variables.}
    \label{fig:5}
  \end{figure}

Therefore, for $m=n$ we consider below $H$ with a vanishingly small value $\ga_m$ and $\ga_1$ significantly larger. Next, to see if the leading $m-1$ singular values of $H$ are distinct, we introduce the probability density function (pdf) of the distribution of singular values of $H$ established in \cite[Th.1]{Shen2001}:
\begin{equation} \label{pdf}
  pdf(\ga_1,\dots,\ga_m)=\frac{1}{C}\exp\left[-\frac{1}{2}\sum_{k=1}^m\ga_k^2\right]
  \prod_{1\leq i<j\leq m}|\ga_j^2-\ga_i^2|,
\end{equation}  
where $C$ normalizes the integral on $[0,\infty)^m.$ It is clear that $pdf(\ga_1,\dots,\ga_m)=0$ whenever $\ga_i=\ga_j$ for some $1\leq i<j\leq m.$ Thus, the singular values of $H$ are almost surely distinct. Indeed, the case $m=n=2$ enables visualization of (\ref{pdf}) in Fig.\ref{fig:5}.

  Apart from depicting the analytical results in (\ref{Pre}) and (\ref{miles}) for $m=n=2$, Fig.\ref{fig:5} shows that the singular values of $H$ are likely to be well-separated. Indeed, for a given size $m=n$ of $H$, one can simply compute an integral of (\ref{pdf}) over a certain integration range to determine the probability of any desired spread of the singular values of~$H.$  

  If we now look at Definition \ref{highill}, the results (\ref{Pre})-(\ref{pdf}) imply that the stochastic linear model (\ref{general}) with a square~$H$ of independent standard normal random variables satisfies Definition \ref{highill} along with the perturbed model (\ref{bubu}) (with a well-defined probability) for $r=m-1$, if:
  \begin{enumerate}
  \item the noise power $\ep$ is such that $\sqrt{\ep}$ is: a) larger than $\ga_m$, b) approximately $\kappa$ times smaller than $\ga_{m-1}$, \label{noi}
  \item one can gauge the difference between $\ga_i$ and $\si_i$ for $i=1,\dots,m-1.$ \label{gau}
  \end{enumerate}
  We assume condition \ref{noi}. below. Regarding \ref{gau}., continuing the work with Gaussian distribution, the norm $\p\De H\p_2$ of a square $\De H$ with independent normal random variables of zero-mean and standard deviation $\si$ is upper-bounded by \cite{Davidson2001}:
  \begin{equation} \label{hd}
    \Pre\left(\p\De H\p_2\leq\si(2\sqrt{m} +t)\right)\geq 1-2\exp(-t^2/2),\ t\geq 0.
  \end{equation}
  For $\si$ several orders of magnitude smaller than 1 and moderate values of $m$, the result (\ref{hd}) ensures that the norm of $\p\De H\p_2$ will be at most one order of magnitude larger than $\si$ with probability almost one [e.g., for $t=10$ in (\ref{hd})]. In such a case, Proposition \ref{YamadaSensei} ensures that the singular values of $H$ and $H+\De H$ will be pairwise separated with a wide margin for $i=1,\dots,m-1.$ This in turn yields values of $\rho_i$ in (\ref{rho}) close to $1$ for $i=1,\dots,m-1$, as promised at the beginning of this section. We will demonstrate the main points of this section in the numerical example in Section \ref{ne}. 

  \section{Numerical Example} \label{ne}
  \subsection{General Setup}
  To illustrate the analysis of Section \ref{cs}, we consider $H\in\mathbb{R}^{4\times 4}$ with independent standard normal random variables such that its least singular value is less than 0.01. According to (\ref{Pre}), this happens in about 1 in 50 cases. As in this paper we consider ill-conditioned settings, we would like $H$ to have a condition number $cond(H)$ of at least $10^3$, which according to (\ref{miles}) happens with non-zero probability (at least 1 in about $2*10^3$ cases) in at most 1 in 50 cases. We note that the probabilities given here refer to a generic matrix with independent standard normal random variables. In certain applications, such as in estimation of brain electrical activity from EEG~/~MEG measurements, the probability of observing highly ill-conditioned array response matrix $H$ is significantly larger \cite{Piotrowski2013}.

  We generated random matrix $H\in\mathbb{R}^{4\times 4}$ possessing the above properties along with perturbation matrix $\De H$, where the entries of $\De H$ are also independent normal variables with zero mean and standard deviation 0.01 such that $\p\De H\p_2=0.034$, \emph{cf.} (\ref{hd}) and discussion below (\ref{hd}). The singular values of $H$ and $H+\De H$ are summarized in Table \ref{tab:svs} below. We consider $\ep=4.928*10^{-4}$ with $\sqrt{\ep}\approx\si_4\approx 0.022.$ From Table \ref{tab:svs} it is seen that this results in signal-to-noise ratio (SNR) value of $SNR[dB]\approx 40$, where SNR is defined as in~(\ref{SNR}). 

  \begin{minipage}{\linewidth}
    \centering
    \bigskip
    \captionof{table}{Singular values of $H$ and $H+\De H$} \label{tab:svs} 
    \begin{tabular}{l | l}\toprule[1.5pt]
      $H$ & $H+\De H$ \\\midrule
      $\ga_1\approx 3.889\approx 175.2\sqrt{\ep}$ &  $\si_1\approx 3.894\approx 175.4\sqrt{\ep}$\\
      $\ga_2\approx 2.426\approx 109.3\sqrt{\ep}$ &  $\si_2\approx 2.435\approx 109.7\sqrt{\ep}$\\
      $\ga_3\approx 0.923\approx 41.6\sqrt{\ep}$ &  $\si_3\approx 0.934\approx 42.1\sqrt{\ep}$\\
      $\ga_4\approx 0.003\approx 0.14\sqrt{\ep}$ &  $\si_4\approx 0.022\approx \sqrt{\ep}$\\
      \bottomrule[1.25pt]
    \end {tabular}\par
    \bigskip
  \end{minipage}

  Hence, such setup satisfies conditions \ref{c1}.-\ref{cend}. of Definition \ref{highill} for $r=3$, and~$\kappa$ in condition \ref{c2}. in Definition \ref{highill} can be set to $\kappa=41.6.$ In particular, it easily satisfies the condition $\kappa=1.75$ of Remark \ref{YS3} in Section \ref{direct}.

  \subsection{Numerical Evaluation of the Results of Section \ref{simple}} \label{ne_simple}
  In this subsection, we set $M=U$ and $N=V$ in (\ref{svdtruH}), \emph{cf.} (\ref{asu}). For $i=1,2,3$, the results of Propositions \ref{wolfram} and \ref{wolfram2} as depicted in Fig.\ref{fig:1}-\ref{fig:2} (up to $10\sqrt{\ep}$) imply that we shall expect a vanishingly small values for each $A_i$ in~(\ref{A_i}). This is indeed the case, as the sum of $A_i$ for $i=1,2,3$ is approximately $1.3*10^{-5}.$ On the other hand, from the results of Propositions~\ref{wolfram} and \ref{wolfram2} as shown in Fig.\ref{fig:3}-\ref{fig:4}, for $i=m=4$ the values of $\ga_4$ and $\si_4$ indicate that we are in the top-left corner of Fig.\ref{fig:4}, which yields value of $B_i$ in (\ref{B_i}) for $i=4$ of $B_4\approx 0.115$, \emph{cf.} Fig.\ref{fig:3}. Indeed, we obtain $\Jn(\wW_{MMSE})\approx 1.1161$ in (\ref{mse_bubu_obtain}), $\Jn(\wW_{r-MMSE})\approx 1.008$ in~(\ref{rrmse_bubu_obtain}), and $\Jn(\wW_{r-SVD})\approx 1.008$ in (\ref{mvp_bubu_obtain}), \emph{cf}. (\ref{R-T})-(\ref{M-T}).

  \subsection{Numerical Evaluation of the Results of Section \ref{ext}}
  We drop the constraint $M=U$ and $N=V$ employed in Subsection \ref{ne_simple}. We obtain the following values of $\rho_i$ (\ref{rho}) in Lemma~\ref{fina} for $i=1,2,3$:
  \begin{equation} \label{cheerio}
    \rho_1\approx 0.9995,\rho_2\approx 0.9995,\rho_3\approx 0.9987.
  \end{equation}
  In particular, we may set $c=0.0013$ in (\ref{ost}) in Remark \ref{R0}. Then, from Theorem \ref{greatinequalities} we obtain in (\ref{elGI}):
  \begin{equation} \label{elGIsim}
    \left\lvert_{(\ref{ext})}\Jn(\wW_{r-MMSE})-_{(\ref{simple})}\Jn(\wW_{r-MMSE})\right\rvert\leq
    0.1404,
  \end{equation}
  in (\ref{terGI}):
  \begin{equation} \label{terGIsim}
    \left\lvert_{(\ref{ext})}\Jn(\wW_{r-SVD})-_{(\ref{simple})}\Jn(\wW_{r-SVD})\right\rvert\leq
    0.1405,
  \end{equation}
  and in (\ref{alGI}):
  \begin{equation} \label{alGIsim}
    \sum_{i=r+1}^m\sum_{j=1}^m(x_i\ga_jk_{i,j}-l_{j,i})^2+\ep\sum_{i=r+1}^mx_i^2=1.019+0.25\geq
    m-r=1.
  \end{equation}
  Thus, from (\ref{elGIsim}) and (\ref{terGIsim}) we have that the MSE values to be achieved without assumption (\ref{asu}) of both $r$-SVD and $r$-MMSE estimators will differ by at most $\approx 0.14$ from their respective MSE values obtained above under assumption~(\ref{asu}). Moreover, in view of (\ref{alGIsim}), from Theorem \ref{greatinequalities} we conclude that the MSE value of the MMSE estimator obtained without assumption (\ref{asu}) will be larger than the corresponding MSE of the $r$-MMSE estimator. Indeed, we obtain $\Jn(\wW_{MMSE})\approx 1.2702$ in~(\ref{sami}), $\Jn(\wW_{r-MMSE})\approx 1.0013$ in (\ref{samija}), and $\Jn(\wW_{r-SVD})\approx 1.0013$ in (\ref{kedira}). 

  Finally, the discussion of Remarks \ref{R1} and \ref{R2} is illustrated numerically in the following tables.

  \begin{minipage}{\linewidth}
    \centering
    \bigskip
    \captionof{table}{Numerical illustration of Remark \ref{R1} on inequality (\ref{elGI})} \label{tab:elGI} 
    \begin{tabular}{l | l}\toprule[1.5pt]
      Remark \ref{R1}, points 1. \& 2. & $2\sum_{i=1}^r\sum_{j=1}^r(1-\rho_i)(x_i\ga_j+1)^2=0.1356$ \\ 
      Remark \ref{R1}, points 1. \& 3. & $2\sum_{i=1}^r\sum_{j=r+1}^m(1-\rho_i)(x_i\ga_j+1)^2=0.0047$ \\ 
      \bottomrule[1.25pt]
    \end {tabular}\par
    \bigskip
  \end{minipage}

  \begin{minipage}{\linewidth}
    \centering
    \bigskip
    \captionof{table}{Numerical illustration of Remark \ref{R1} on inequality (\ref{terGI})} \label{tab:terGI} 
    \begin{tabular}{l | l}\toprule[1.5pt]
      Remark \ref{R1}, points 1. \& 2. & $2\sum_{i=1}^r\sum_{j=1}^r(1-\rho_i)(\si_i^{-1}\ga_j+1)^2=0.1357$ \\ 
      Remark \ref{R1}, points 1. \& 3. & $2\sum_{i=1}^r\sum_{j=r+1}^m(1-\rho_i)(\si_i^{-1}\ga_j+1)^2=0.0047$ \\ 
      \bottomrule[1.25pt]
    \end {tabular}\par
    \bigskip
  \end{minipage}

  \begin{minipage}{\linewidth}
    \centering
    \bigskip
    \captionof{table}{Numerical illustration of Remark \ref{R2} on inequality (\ref{alGI})} \label{tab:alGI} 
    \begin{tabular}{l | l}\toprule[1.5pt]
      Remark \ref{R2}, point 1. & $\sum_{i=r+1}^m\sum_{j=1}^r(x_i\ga_jk_{i,j}-l_{j,i})^2=0.1537$ \\ 
      Remark \ref{R2}, point 2. & $\sum_{i=r+1}^m\sum_{j=r+1}^m(x_i\ga_jk_{i,j}-l_{j,i})^2=0.8652$ \\ 
      Remark \ref{R2}, point 3. & $\ep\sum_{i=r+1}^mx_i^2=0.25$ \\ 
      \bottomrule[1.25pt]
    \end {tabular}\par
    \bigskip
  \end{minipage}
  
  \section{Conclusion} \label{conclusion}
  By perturbing the array response matrix, we derived in Sections \ref{simple} and~\ref{ext} explicit sufficient conditions under which reduced-rank estimators achieve lower MSE than the theoretically MSE-optimal MMSE estimator. The main findings of the paper were verified by numerical simulations. The future research will focus on applications of the derived results in signal processing applications, and extending the proposed approach to ridge regression and shrinkage methods. Indeed, the simplest ridge regression estimator may be expressed in our settings as [\emph{cf}. (\ref{mmseish})-(\ref{mvpform_s})]:
  \begin{equation} \label{RR}
    \widetilde{W}_{RR}=(H^tH+\eta I_m)^{-1}H^t=V(\Si^t\Si+\eta I_m)^{-1}\Si^tU^t,\ \eta>0.
  \end{equation}
  Then, one can proceed with establishing the MSE expressions of (\ref{RR}) for both simplified and generic perturbations as has been done in Sections \ref{simple} and~\ref{ext}, respectively, for the MMSE, reduced-rank MMSE, and truncated SVD estimators.

  \appendix

  \section{Proof of Proposition \ref{YamadaSensei}} \label{pd_YamadaSensei}
  With notation as in Fact \ref{Stewart} in \ref{kru}, and for $A=H$, $E=\De H$, and $\widetilde{A}=H+\De H$, one has from (\ref{sv_bounds}) that:
  \begin{equation} \label{sensei}
    \si_i^2=(\ga_i+\up_i)^2+\eta_i^2,\ i=1,\dots,m.
  \end{equation}
  We use the fact that the singular values of orthogonal projection matrices are either 0 or 1. Then, using (\ref{ep_i}) one has that:
  \begin{equation} \label{s1}
    |\up_i|\leq\p P_{\ra{H}}\De H\p_2\leq\p\De H\p_2,\ i=1,\dots,m,
  \end{equation}
  and similarly, using (\ref{eta_i}) one has that:
  \begin{equation} \label{s2}
    0\leq\ds{\min}_2(P_{\ra{H}}^\perp\De H)\leq\eta_i\leq\p P_{\ra{H}}^\perp \De H\p_2\leq\p\De H\p_2,\ i=1,\dots,m.
  \end{equation}
  Thus, under assumption (\ref{YSasu}), from (\ref{sensei}) we obtain that:
  \begin{equation} \label{kr}
    \si_i\geq\sqrt{(\ga_i-\p\De H\p_2)^2+0}=\ga_i-\p\De H\p_2,
  \end{equation}
  and
  \begin{multline} \label{sbe}
    \si_i\leq\sqrt{(\ga_i+\p\De H\p_2)^2+\p\De H\p_2^2}\leq\sqrt{(\ga_i+\sqrt{2}\p\De H\p_2)^2}=\\\ga_i+\sqrt{2}\p\De H\p_2.\
  \end{multline}

  \section{Proof of Proposition \ref{s_expr}} \label{pd_s_expr}
  From (\ref{mse}) and (\ref{mmseish}) we have:
  \begin{multline*}
    \Jn(\wW_{MMSE})=tr[\wW_{MMSE}\Ry (\wW_{MMSE})^t]-2tr[\wW_{MMSE}\Ryx]+tr[\Rx]=\\
    tr[V\Si^t(\Si\Si^t+\ep I_n)^{-1}U^t*U(\Ga\Ga^t+\ep I_n)U^t*U(\Si\Si^t+\ep I_n)^{-1}\Si V^t]-\\2tr[V\Si^t(\Si\Si^t+\ep I_n)^{-1}U^t*U\Ga V^t]+tr[I_m]=\\\sum_{i=1}^m{{\si_i^2(\ga_i^2+\ep)}\over
      {(\si_i^2+\ep)^2}}-
    2\sum_{i=1}^m{{\ga_i\si_i}\over{\si_i^2+\ep}}+m=\sum_{i=1}^mx_i^2(\ga_i^2+\ep)-
    2\sum_{i=1}^mx_i\ga_i+m.
  \end{multline*}
  Similarly, from (\ref{mse}) and (\ref{rrmmseish}) we have:
  \begin{multline*}
    \Jn(\wW_{r-MMSE})=\\tr[\wW_{r-MMSE}\Ry (\wW_{r-MMSE})^t]-2tr[\wW_{r-MMSE}\Ryx]+tr[\Rx]=\\
    tr[VI_m^r\Si^t(\Si\Si^t+\ep I_n)^{-1}U^t*U(\Ga\Ga^t+\ep I_n)U^t*U(\Si\Si^t+\ep I_n)^{-1}\Si I_m^rV^t]-\\2tr[VI_m^r\Si^t(\Si\Si^t+\ep I_n)^{-1}U^t*U\Ga V^t]+tr[I_m]=\\\sum_{i=1}^r{{\si_i^2(\ga_i^2+\ep)}\over
      {(\si_i^2+\ep)^2}}-
    2\sum_{i=1}^r{{\ga_i\si_i}\over{\si_i^2+\ep}}+m=\sum_{i=1}^rx_i^2(\ga_i^2+\ep)-
    2\sum_{i=1}^rx_i\ga_i+m,
  \end{multline*}
  and from (\ref{mse}) and (\ref{mvpform_s}) we have:
  \begin{multline*}
    \Jn(\wW_{r-SVD})=tr[\wW_{r-SVD}\Ry (\wW_{r-SVD})^t]-2tr[\wW_{r-SVD}\Ryx]+tr[\Rx]=\\
    tr\left[V\Si_r^\dagger U^t*U(\Ga\Ga^t+\ep I_n)U^t*U(\Si_r^\dagger)^t V^t\right]-\\2tr\left[V\Si_r^\dagger U^t*U\Ga V^t\right]+tr[I_m]=\\\sum_{i=1}^r{{\ga_i^2+\ep}\over
      {\si_i^2}}-2\sum_{i=1}^r{{\ga_i}\over{\si_i}}+m=m-r+\ep\sum_{i=1}^r{{1}\over{\si_i^2}}+\sum_{i=1}^r\biggl({{\ga_i}\over{\si_i}}-1\biggr)^2.\
  \end{multline*}

  \section{Proof of Proposition \ref{elilama}} \label{pd_elilama}
  With $x_i=\frac{\si_i}{\si_i^2+\ep}$ defined in (\ref{ixi}) one has:
  \begin{equation} \label{terribly_partial}
    \frac{\partial A_i}{\partial\ga_i}=2\ga_i(x_i^2-\si_i^{-2})
    -2(x_i-\si_i^{-1}),\ i=1,\dots,m, 
  \end{equation}
  and
  \begin{equation} \label{He}
    \frac{\partial^2 A_i}{\partial\ga_i^2}=2(x_i^2-\si_i^{-2})<0,\ i=1,\dots,m. 
  \end{equation}
  Similarly:
  \begin{equation} \label{terribly_partial_B}
    \frac{\partial B_i}{\partial\ga_i}=2\ga_ix_i^2
    -2x_i,\ i=1,\dots,m,
  \end{equation}
  and
  \begin{equation}
    \frac{\partial^2 B_i}{\partial\ga_i^2}=2x_i^2>0,\ i=1,\dots,m.
  \end{equation}
  By equating (\ref{terribly_partial}) to 0 it is straightforward to verify that the maximum of $A_i$ is achieved~at
  \begin{equation} 
    \ga_i^{A_i^{\min}}={{\ep\si_i(\si_i^2+\ep)}\over{(\si_i^2+\ep)^2-\si_i^4}}=
       {{\si_i(\si_i^2+\ep)}\over{2\si_i^2+\ep}},\ i=1,\dots,m,
  \end{equation}
  while equating (\ref{terribly_partial_B}) to 0 reveals that the minimum of $B_i$ is achieved at 
  \begin{equation} 
    \ga_i^{B_i^{\min}}=\si_i+{{\ep}\over{\si_i}}=x_i^{-1},\ i=1,\dots,m.\
  \end{equation}

  \section{Proof of Proposition \ref{wolfram}} \label{pd_wolfram}
  Due to concavity of $A_i$ and convexity of $B_i$, the proof reduces to finding zero-crossings of $A_i$ and $B_i$ as functions of $\ga_i.$ To preserve space, this task has been completed using Wolfram Mathematica 9 (Mathematica) \emph{Reduce} function, followed by \emph{FullSimplify} function to achieve simplest algebraical expressions in (\ref{whew1eq}) and (\ref{whew2eq}).

  \section{Proof of Proposition \ref{wolfram2}} \label{pd_wolfram2}
  The fact that parametrization $\ga_i=a_{\ga_i}\sqrt{\ep}$ and $\si_i=a_{\si_i}\sqrt{\ep}$ allows to eliminate $\ep$ from (\ref{A_i}) and (\ref{B_i}) is obvious. The expressions (\ref{ja}) and (\ref{pan}) are obtained from (\ref{A_i}) and (\ref{B_i}), respectively, by simple but tedious algebraic manipulations of the expressions of $A_i$ and $B_i.$ To preserve space, this task has been completed using Mathematica.

  We prove now that (\ref{kamieni}) is the unique maximum of $A_i$ for fixed $\ep.$ From Proposition \ref{elilama} we obtain that $\ga_i$ must be of the form (\ref{ingen}) if maximum of $A_i$ is to be achieved. Algebraic manipulations show that $A_i$ can be expressed for such argument as
  \begin{equation} \label{Amax1}
    A_i\left({{\si_i(\si_i^2+\ep)}\over{2\si_i^2+\ep}},\si_i,\ep\right)=
    \ep\left(x_i^2+\frac{1}{\ep+2\si_i^2}-\frac{1}{\si_i^2}\right),
  \end{equation}
  where $x_i=\frac{\si_i}{\si_i^2+\ep}$ is defined in (\ref{ixi}). Introducing parametrization $\si_i=a_{\si_i}\sqrt{\ep}$ for $a_{\si_i}>0$ allows to eliminate $\ep$ from (\ref{Amax1}), leaving $A_i$ evaluated for $\ga_i$ of the form (\ref{ingen}) only as a function of $a_{\si_i}.$ Namely, after some algebraic manipulations we find that (\ref{Amax1}) can be expressed as:
  \begin{equation} \label{Amax2}
    A_i\left({{a_{\si_i}(a_{\si_i}^2+1)}\over
      {2a_{\si_i}^2+1}}\sqrt{\ep},a_{\si_i}\sqrt{\ep},\ep\right)=
    \frac{-1-3a_{\si_i}^2-2a_{\si_i}^4+a_{\si_i}^6}
         {(1+2a_{\si_i}^2)(a_{\si_i}+a_{\si_i}^3)^2}.
  \end{equation}
  Calculating the derivative of (\ref{Amax2}) shows that for positive $a_{\si_i}$ it vanishes only for $a_{\si_i}\approx 2.611$, and the second derivative of (\ref{Amax2}) is negative at this point. Thus, for fixed $\ep$, $A_i$ has the unique maximum for $\si_i=2.611\sqrt{\ep}$ and 
  \begin{equation*}
    \ga_i={{a_{\si_i}(a_{\si_i}^2+1)}\over
      {2a_{\si_i}^2+1}}\sqrt{\ep}=1.394\sqrt{\ep},
  \end{equation*}
  with $A_i(1.394\sqrt{\ep},2.611\sqrt{\ep},\ep)\approx 0.033.$ 

  Let us now restrict the domain of $\ga_i$ to the interval $\left(0,\ x_i^{-1}-\be_i^B\right).$ In view of Remark~\ref{YS2}, any maximum of~$B_i$ as a function of $\ga_i$ and $\si_i$ may be achieved only for $\ga_i\to 0.$ We have:
  \begin{equation} \label{Bmax1}
    \lim_{\ga_i\rightarrow 0}B_i(\ga_i,\si_i,\ep)=\ep x_i^2,
  \end{equation}
  which can be expressed as a function of $a_{\si_i}$ only:
  \begin{equation} \label{Bmax2}
    \lim_{\ga_i\rightarrow 0}B_i(\ga_i,a_{\si_i}\sqrt{\ep},\ep)=\frac{a_{\si_i}^2}{(1+a_{\si_i}^2)^2}. \end{equation}
  Calculating the derivative of (\ref{Bmax2}) shows that for nonnegative $a_{\si_i}$ it vanishes for $a_{\si_i}=1$ and $a_{\si_i}=0$, with the second derivative of (\ref{Bmax2}) negative for $a_{\si_i}=1$ and positive for $a_{\si_i}=0.$ Thus, for fixed $\ep$ and $\si_i=\sqrt{\ep}$ we have
  \begin{equation} \label{kupaconfirmed}
    \sup_{\ga_i\in\left(0,\ x_i^{-1}-\be_i^B\right)}B_i(\ga_i,\sqrt{\ep},\ep)=\lim_{\ga_i\rightarrow 0}B_i(\ga_i,\sqrt{\ep},\ep)=0.25.\
  \end{equation}

  \section{Proof of Proposition \ref{mse_exprs}} \label{pd_mse_exprs}
  With notation introduced in (\ref{no1})-(\ref{L}), we first insert (\ref{mmseish}) expressed as $\wW_{MMSE}=VXU^t$ into (\ref{mse}):
  \begin{multline} 
    \Jn(\wW_{MMSE})=tr[VXU^tMYM^tUX^tV^t]-
    2tr[VXU^tM\Ga N^t]+m=\\tr[X^tXKYK^t]-2tr[XK\Ga L]+m.
  \end{multline}
  The diagonal entries of $KYK^t\in\nn$ are of the form:
  \begin{multline}
    (KYK^t)_{i,i}=\sum_{j=1}^mk_{i,j}(\ga_j^2+\ep)k_{i,j}+\sum_{j={m+1}}^nk_{i,j}\ep
    k_{i,j}=\\\sum_{j=1}^mk_{i,j}^2\ga_j^2+\ep\sum_{j=1}^nk_{i,j}^2=\sum_{j=1}^mk_{i,j}^2\ga_j^2+\ep,\  i=1,\dots,n,
  \end{multline}
  since $\sum_{j=1}^nk_{i,j}^2=1.$ Hence, with (\ref{ixi}), the diagonal entries of $X^tXKYK^t\in\nn$ are of the form:
  \begin{equation} \label{profi1}
    (X^tXKYK^t)_{i,i}=\left\{
    \begin{array}{ll}
      x_i^2\left(\sum_{j=1}^mk^2_{i,j}\ga^2_j+\ep\right) & 1\leq i\leq m\vspace{0.1cm}\\
      0 & m<i\leq n,\\
    \end{array}\right.
  \end{equation}
  and we obtain therefore that  $tr[X^tXKYK^t]=\sum_{i=1}^m\left[x_i^2\left(\sum_{j=1}^mk_{i,j}^2\ga_j^2+\ep\right)\right].$ Similarly, the elements on the main diagonal of $K\Ga L\in\nm$ are of the form:
  \begin{equation}
    (K\Ga L)_{i,i}=\sum_{j=1}^mk_{i,j}\ga_jl_{j,i}=\sum_{j=1}^m\ga_jk_{i,j}l_{j,i},\  i=1,\dots,m,
  \end{equation}
  and thus diagonal entries of $XK\Ga L\in\mm$ are of the form:
  \begin{equation} \label{profi2}
    (XK\Ga L)_{i,i}=x_i\sum_{j=1}^m\ga_jk_{i,j}l_{j,i},\  i=1,\dots,m,
  \end{equation}
  therefore $tr[XK\Ga L]=\sum_{i=1}^m\left(x_i\sum_{j=1}^m\ga_jk_{i,j}l_{j,i}\right)$, which completes the proof of (\ref{sami}). 

  The proof of (\ref{samija}) is obtained as follows: from (\ref{rrmmseish}) we have that $\wW_{r-MMSE}=VI^r_mXU^t.$ Inserting it into (\ref{mse}) shows that
  \begin{multline} \label{bread}
    \Jn(\wW_{r-MMSE})=tr[VI^r_mXU^tMYM^tUX^tI^r_mV^t]-
    2tr[VI^r_mXU^tM\Ga N^t]+m=\\ tr[X^tI^r_mXKYK^t]-2tr[I^r_mXK\Ga L]+m.
  \end{multline}
  The diagonal matrix $X^tI^r_mX\in\nn$ has diagonal entries $x_i^2$ for $i=1,\dots,r$ and $0$ for $i=r+1,\dots,n.$ Similarly, the main diagonal of $I^r_mX\in\mn$ is the vector $(x_1,\dots,x_r,0,\dots,0)\in\sm.$ Hence, the proof of (\ref{samija}) follows immediately from (\ref{profi1}) and (\ref{profi2}).

  The proof of (\ref{kedira}) follows by inserting (\ref{mvpform_s}) into (\ref{mse}):
  \begin{multline} \label{fan}
    \Jn(\wW_{r-SVD})=
    tr[V\Si^\dagger_rU^tMYM^tU(\Si^\dagger_r)^tV^t]-2tr[V\Si^\dagger_rU^tM\Ga N^t]+m=\\ tr[\Si^\dagger_rKYK^t(\Si^\dagger_r)^t]-2tr[\Si^\dagger_rK\Ga L]+m=\\tr[(\Si^\dagger)^tI^r_m\Si^\dagger KYK^t]-2tr[I^r_m\Si^\dagger K\Ga L]+m,
  \end{multline}
  where $\Si^\dagger\in\mn$ contains on its main diagonal the reciprocals of the singular values of $H$ and zeros elsewhere. Note that (\ref{fan}) has the same form as (\ref{bread}), with $X$ replaced by $\Si^\dagger$ in the expression on the rightmost side in (\ref{fan}). Thus, the proof of (\ref{kedira}) follows exactly the same lines as the proof of (\ref{samija}).

  \section{Proof of Lemma \ref{fina}} \label{pd_fina}
  The proof is divided into three parts.
  \paragraph{Expressions for $k_{i,i}$ and $l_{i,i}$} We recall first that the columns of orthogonal matrices are of unit norm. Then, from (\ref{K}) and (\ref{L}) we obtain, respectively, that
  \begin{equation} \label{cos1}
    k_{i,i}=u_i^tm_i=m_i^tu_i=\frac{m_i^tu_i}{\p m_i\p\ \p u_i\p}=\cos\angle(m_i,u_i),\  i\in\{1,\dots,m\},
  \end{equation}
  and
  \begin{equation} \label{cos2}
    l_{i,i}=n_i^tv_i=\frac{n_i^tv_i}{\p n_i\p\ \p v_i\p}=\cos\angle(n_i,v_i),\  i\in\{1,\dots,m\},
  \end{equation}
  thus $2\geq\lvert k_{i,i}\rvert+\lvert l_{i,i}\rvert\geq k_{i,i}^2+l_{i,i}^2.$ Therefore, to prove (\ref{faaaa}) it suffices to show that $k_{i,i}^2+l_{i,i}^2\geq 2\rho_i.$
  \paragraph{Inequality $k_{i,i}^2+l_{i,i}^2\geq 2\rho_i$} Let us first take a look into conditions \ref{c1}.-\ref{cend}. of Definition \ref{highill}. Then, for $A=H$, $E=\De H$, and $\widetilde{A}=H+\De H$, from Definition~\ref{master_def} and Proposition \ref{YamadaSensei} [considered via condition \ref{c1}. in Definition \ref{highill}] it is seen that $\de_{i,i}$ in (\ref{Wedin_delta}) in Wedin's Theorem in Fact \ref{Wedin} in \ref{kru} is equal to $\de_i$ defined through (\ref{mex})-(\ref{08}) for $i=1,\dots,r-1.$ For $i=r$, the form of $\bs{\ga_{r,NEXT}}$ in (\ref{08}) for $\de_r$ in (\ref{mex}) is obtained by combining conditions \ref{c1}.~-~\ref{cend}. of Definition \ref{highill}. Therefore, with $R_i=(H-(H+\De H))v_i=-(\De H)v_i$, $S_i=(H^t-(H+\De H)^t)u_i=-(\De H)^tu_i$, one obtains from Wedin's Theorem for $j=i$ that
  \begin{multline} \label{dzizas}
    \sin^2\angle(m_i,u_i)+\sin^2\angle(n_i,v_i)\leq{{\p R_i\p^2+\p S_i\p^2}\over{\de_i^2}}\leq\\{{\p\De H\p_2^2}(\p v_i\p^2+\p u_i\p^2)\over{\de_i^2}}={{2\p\De H\p_2^2}\over{\de_i^2}},\ i=1,\dots,r,
  \end{multline}
  where $\de_i$ is defined through (\ref{mex})-(\ref{08}). The last inequality in (\ref{faaaa}) follows now from (\ref{dzizas}) using Pythagorean trigonometric identity. 
  \paragraph{Inequalities $1\geq k_{i,i}l_{i,i}>0$} The fact that $1\geq k_{i,i}l_{i,i}$ is obvious from (\ref{cos1}) and (\ref{cos2}). Now, for $A=H$, $E=\De H$, and $\widetilde{A}=H+\De H$, from Dopico's Theorem as given in Fact \ref{Dopico} in \ref{kru} one has for $j=i$, $\zeta_i=\zeta_{i,i}$ and $i=1,\dots,r$ that:
  \begin{multline} \label{fia}
    \min_{x\in\{-1,1\}}\left\{\p xm_i-u_i\p^2+\p xn_i-v_i\p^2\right\}\leq
    2{{\p R_i\p^2+\p S_i\p^2}\over{\zeta_i^2}}\leq\\ 2{{\p\De H\p_2^2}(\p v_i\p^2+\p u_i\p^2)\over{\zeta_i^2}}={{4\p\De H\p_2^2}\over{\zeta_i^2}},
  \end{multline}
  where $\zeta_i=\min\{\de_i,\si_i+\ga_i\}.$ However, from condition \ref{c1}. in Definition \ref{highill} it is seen that $\si_i>\ga_{i+1}$ for $i=1,\dots,r-1$, and by combining conditions \ref{c1}.-\ref{cend}. in Definition \ref{highill} it is seen that also $\si_r>\ga_{r+1}.$ Therefore, we obtain that $\de_i$ in (\ref{mex}) is such that $\de_i<\si_i$ and consequently $\zeta_i=\de_i$ in (\ref{fia}) for $i=1,\dots,r.$ We consider now two cases in (\ref{fia}).
  \begin{itemize}
  \item $\displaystyle{\min_{x\in\{-1,1\}}}\left\{\p xm_i-u_i\p^2+\p xn_i-v_i\p^2\right\}=\p m_i-u_i\p^2+\p n_i-v_i\p^2.$ Then, in view of the above one has:
    \begin{multline}
      \p m_i-u_i\p^2+\p n_i-v_i\p^2=\\\p m_i\p^2+\p u_i\p^2-2k_{i,i}+\p n_i\p^2+\p v_i\p^2-2l_{i,i}=\\4-2(k_{i,i}+l_{i,i})\leq {{4\p\De H\p_2^2}\over{\delta_i^2}},
    \end{multline}
    thus
    \begin{equation}
      k_{i,i}+l_{i,i}\geq 2\left(1-{{\p\De H\p_2^2}\over{\delta_i^2}}\right)=2\rho_i.
    \end{equation}
    Therefore, if $\rho_i>1/2$, then $k_{i,i}+l_{i,i}>1$, hence it must be $k_{i,i}>0$ and $l_{i,i}>0$, and consequently $k_{i,i}l_{i,i}>0.$
  \item $\displaystyle{\min_{x\in\{-1,1\}}}\left\{\p xm_i-u_i\p^2+\p xn_i-v_i\p^2\right\}=\p -m_i-u_i\p^2+\p -n_i-v_i\p^2=\p m_i+u_i\p^2+\p n_i+v_i\p^2.$ Similarly as above, one has in such a case that:
    \begin{multline}
      \p m_i+u_i\p^2+\p n_i+v_i\p^2=\\\p m_i\p^2+\p u_i\p^2+2k_{i,i}+\p n_i\p^2+\p v_i\p^2+2l_{i,i}=\\4+2(k_{i,i}+l_{i,i})\leq {{4\p\De H\p_2^2}\over{\delta_i^2}},
    \end{multline}
    thus
    \begin{equation}
      k_{i,i}+l_{i,i}\leq -2\left(1-{{\p\De H\p_2^2}\over{\delta_i^2}}\right)=-2\rho_i.
    \end{equation}
    Therefore, if $\rho_i>1/2$, then $k_{i,i}+l_{i,i}<-1$, hence it must be $k_{i,i}<0$ and $l_{i,i}<0$, and consequently $k_{i,i}l_{i,i}>0.$ This completes the proof of (\ref{fuuuu}) for $\rho_i>1/2$ for $i\in\{1,\dots,r\}.$
  \end{itemize}

  \section{Proof of Proposition \ref{funda}} \label{pd_funda}
  The first inequality in (\ref{conc0}) follows from $k_{i,i}^2+l_{i,i}^2\geq 2\rho_i$ in (\ref{faaaa}), since $k_{i,i}^2\geq 2\rho_i-l_{i,i}^2\geq 2\rho_i-1$, in view of $l_{i,i}^2\leq 1.$ The second inequality in (\ref{conc0}) can be proved analogously. The third inequality in (\ref{conc0}) follows from the first two, since $|k_{i,i}|\geq\sqrt{2\rho_i-1}$ and $|l_{i,i}|\geq\sqrt{2\rho_i-1}$, thus it must be $|k_{i,i}l_{i,i}|\geq 2\rho_i-1$, and from (\ref{fuuuu}) we obtain that $|k_{i,i}l_{i,i}|=k_{i,i}l_{i,i}$ for $\rho_i>1/2.$

  The first inequality in (\ref{conc1}) can be proved as follows, taking into account that $\sum_{j=1}^nk_{i,j}^2=1$:
  \begin{equation}
    \sum_{j=1,j\neq i}^nk_{i,j}^2=1-k_{i,i}^2\leq 1-(2\rho_i-1)=2(1-\rho_i).
  \end{equation}
  The second inequality in (\ref{conc1}) can be proved analogously. The third and fourth inequalities in (\ref{conc1}) are now obvious. One also has $\sum_{j=1,j\neq i}^mk_{i,j}^2\leq 2(1-\rho_i)$ and $\sum_{j=1,j\neq i}^mk_{j,i}^2\leq 2(1-\rho_i)$, since $m\leq n.$ Furthermore, the same arguments as above can be used to prove inequalities in (\ref{conc1.5}).

  The first inequality in (\ref{conc2}) follows from Cauchy-Schwarz inequality:
  \begin{equation}
    \left|\sum_{j=1,j\neq i}^mk_{i,j}l_{j,i}\right|\leq\sqrt{\sum_{j=1,j\neq i}^mk_{i,j}^2}\sqrt{\sum_{j=1,j\neq i}^ml_{j,i}^2}\leq\sqrt{2(1-\rho_i)}\sqrt{2(1-\rho_i)}=2(1-\rho_i).
  \end{equation}
  The second inequality in (\ref{conc2}) is obvious in view of $|k_{i,j}l_{j,i}|=|k_{i,j}||l_{j,i}|$, and the fact that $k_{i,j}^2\leq 2(1-\rho_i)$ and $l_{j,i}^2\leq 2(1-\rho_i).$

  \section{Proof of Theorem \ref{greatinequalities}} \label{pd_greatinequalities}
  In view of (\ref{el}) and (\ref{ter}) in Corollary \ref{C_main}, it is sufficient to prove (\ref{elGI}), as the proof of (\ref{terGI}) will follow by replacing $x_i$ with $\si_i^{-1}$ in the proof of (\ref{elGI}). From (\ref{phipsi}), (\ref{el}) and the fact that $\sum_{j=1}^ml_{j,i}^2=1$ for $i=1,\dots,r$ one has:
  \begin{multline} \label{Gen0}
    \left\lvert_{(\ref{ext})}\Jn(\wW_{r-MMSE})-_{(\ref{simple})}\Jn(\wW_{r-MMSE})\right\rvert=\\\left\lvert
    \sum_{i=1}^r\left[x_i^2\left(\phi_i^2-\ga_i^2\right)+
      2x_i\left(\ga_i-\psi_i\right)\right]\right\rvert=\\\left\lvert
    \sum_{i=1}^r\left[x_i^2\left(\sum_{j=1}^mk_{i,j}^2\ga_j^2-\ga_i^2\right)+
      2x_i\left(\ga_i-\sum_{j=1}^m\ga_jk_{i,j}l_{j,i}\right)\right]\right\rvert=\\
    \left\lvert
    \sum_{i=1}^r\left[\sum_{j=1}^m(x_i\ga_jk_{i,j}-l_{j,i})^2-1-(x_i\ga_i-1)^2+1\right]\right\rvert
    \leq\\\sum_{i=1}^r\left\lvert\sum_{j=1}^m(x_i\ga_jk_{i,j}-l_{j,i})^2-(x_i\ga_i-1)^2\right\rvert
    \leq\\\sum_{i=1}^r\sum_{j=1,j\neq i}^m(x_i\ga_jk_{i,j}-l_{j,i})^2+
    \sum_{i=1}^r\left\lvert(x_i\ga_ik_{i,i}-l_{i,i})^2-(x_i\ga_i-1)^2\right\rvert.
  \end{multline}
  From Proposition \ref{funda} we obtain in particular that $k_{i,j}^2\leq 2(1-\rho_i)$ and $l_{j,i}^2\leq 2(1-\rho_i)$ for $i=1,\dots,r$, $j=1,\dots,m$ and $j\neq i.$ Therefore
  \begin{multline} \label{Gen1}
    \sum_{i=1}^r\sum_{j=1,j\neq i}^m(x_i\ga_jk_{i,j}-l_{j,i})^2\leq
    \sum_{i=1}^r\sum_{j=1,j\neq i}^m\left(x_i\ga_j\sqrt{2(1-\rho_i)}+\sqrt{2(1-\rho_i)}\right)^2=\\
    2\sum_{i=1}^r\sum_{j=1,j\neq i}^m(1-\rho_i)(x_i\ga_j+1)^2.
  \end{multline}
  Similarly, from Proposition \ref{funda} we also have that $k_{i,i}^2\geq 2\rho_i-1$, $l_{i,i}^2\geq 2\rho_i-1,$ and  $k_{i,i}l_{i,i}\geq 2\rho_i-1$ for $i=1,\dots,r$, therefore
  \begin{multline} \label{Gen2}
    \sum_{i=1}^r\left\lvert(x_i\ga_ik_{i,i}-l_{i,i})^2-(x_i\ga_i-1)^2\right\rvert=\\
    \sum_{i=1}^r\left\lvert x_i^2\ga_i^2k_{i,i}^2-2x_i\ga_ik_{i,i}l_{i,i}+l_{i,i}^2-x_i^2\ga_i^2+2x_i\ga_i-1\right\rvert=\\
    \sum_{i=1}^r\left\lvert x_i^2\ga_i^2(k_{i,i}^2-1)+2x_i\ga_i(1-k_{i,i}l_{i,i})+l_{i,i}^2-1\right\rvert\leq\\
    \sum_{i=1}^r\left[x_i^2\ga_i^2\lvert k_{i,i}^2-1\rvert+2x_i\ga_i\lvert 1-k_{i,i}l_{i,i}\rvert+\lvert l_{i,i}^2-1\rvert\right]=\\
    \sum_{i=1}^r\left[x_i^2\ga_i^2\lvert 1-k_{i,i}^2\rvert+2x_i\ga_i\lvert 1-k_{i,i}l_{i,i}\rvert+\lvert 1-l_{i,i}^2\rvert\right]\leq\\
    \sum_{i=1}^r\left[x_i^2\ga_i^2(1-2\rho_i+1)+2x_i\ga_i(1-2\rho_i+1)+(1-2\rho_i+1)\right]=\\
    2\sum_{i=1}^r(1-\rho_i)(x_i^2\ga_i^2+2x_i\ga_i+1)=2\sum_{i=1}^r(1-\rho_i)(x_i\ga_i+1)^2.
  \end{multline}
  By combining (\ref{Gen1}) and (\ref{Gen2}), in view of (\ref{Gen0}) the proof of (\ref{elGI}) is completed. 

  Consider now (\ref{al}) in Corollary \ref{C_main}. Similarly as above one has:
  \begin{multline} \label{Gen3}
    _{(\ref{ext})}\Jn(\wW_{MMSE})-_{(\ref{ext})}\Jn(\wW_{r-MMSE})=
    \sum_{i=r+1}^mx_i^2(\phi_i^2+\ep)-2\sum_{i=r+1}^mx_i\psi_i=\\
    \sum_{i=r+1}^m(x_i^2\phi_i^2-2x_i\psi_i)+\ep\sum_{i=r+1}^mx_i^2=\\
    \sum_{i=r+1}^m\left(x_i^2\sum_{j=1}^mk_{i,j}^2\ga_j^2-
    2x_i\sum_{j=1}^m\ga_jk_{i,j}l_{j,i}\right)+\ep\sum_{i=r+1}^mx_i^2=\\
    \sum_{i=r+1}^m\left(\sum_{j=1}^m(x_i\ga_jk_{i,j}-l_{j,i})^2-1\right)+\ep\sum_{i=r+1}^mx_i^2=\\
    \sum_{i=r+1}^m\sum_{j=1}^m(x_i\ga_jk_{i,j}-l_{j,i})^2-(m-r)+\ep\sum_{i=r+1}^mx_i^2,
  \end{multline}
  from which condition (\ref{alGI}) follows.

  \section{Known Results Used} \label{kru}
  \begin{fact}[\cite{Stewart1979}] \label{Stewart}
    Let matrix $A\in\nm$ have singular values $\bs{\alpha}=(\alpha_1,\dots,\alpha_m)$, and matrix $\widetilde{A}=A+E$ for some $E\in\nm$ have singular values $\bs{\widetilde{\alpha}}=(\tilde{\alpha}_1,\dots,\tilde{\alpha}_m)$, organized in both cases in nonincreasing order. Then:
    \begin{equation} \label{sv_bounds}
      \tilde{\alpha}_i^2=(\alpha_i+\up_i)^2+\eta_i^2,\ i=1,\dots,m,
    \end{equation}
    where
    \begin{equation} \label{ep_i}
      |\up_i|\leq\p P_{\ra{A}}E\p_2,
    \end{equation}
    and
    \begin{equation} \label{eta_i}
      \ds{\min}_2(P_{\ra{A}}^\perp E)\leq\eta_i\leq\p P_{\ra{A}}^\perp E\p_2,
    \end{equation}
    where  $\p\bullet\p_2$ is the matrix spectral norm (the largest singular value of a matrix) \cite{Horn1985}, $P_{\ra{A}}$ is the orthogonal projection matrix onto the range of $A$, $P_{\ra{A}}^\perp$ is the orthogonal projection matrix onto the orthogonal complement of the range of $A$, and $\ds{\min}_2(P_{\ra{A}}^\perp E)$ is the smallest singular value of $P_{\ra{A}}^\perp E.$
  \end{fact}

  \begin{fact}[$\mbox{Wedin's Theorem \cite[p.260]{Stewart1990}}$] \label{Wedin}
    Let $A$ and $\widetilde{A}=A+E$ for some $E\in\nm$ be two $n\times m$ matrices such that $n\geq m$, both with distinct singular values. For $i,j\in\{1,\dots,m\}$, let us pick $\alpha_i=i$-th singular value of $A$ along with $m_i$ and $n_i$: left and right singular vectors corresponding to $\alpha_i$, and similarly $\tilde{\alpha}_j=j$-th singular value of $\tilde{A}$ along with $u_j$ and $v_j$: left and right singular vectors corresponding to $\tilde{\alpha}_j.$ Let us denote the set of singular values of~$A$ as  $\bs{\alpha_{SET}}=\{\alpha_1,\dots,\alpha_m\}.$ Furthermore, let us define:
    \begin{equation} \label{Wedin_delta}
      \de_{i,j}=\min_{{\alpha}\in\bs{\alpha_{i, SET-}}}|\tilde{\alpha}_j-\alpha|,
    \end{equation}
    where $\bs{\alpha_{i, SET-}}=\bs{\alpha_{SET}}\setminus\{\alpha_i\}\cup\{0\}$ if $n>m$ and $\bs{\alpha_{i, SET-}}=\bs{\alpha_{SET}}\setminus\{\alpha_i\}$ if $n=m.$ Moreover, let us set:
    \begin{equation*}
      R_j=(A-\widetilde{A})v_j
    \end{equation*}
    and 
    \begin{equation*}
      S_j=(A^t-\widetilde{A}^t)u_j.
    \end{equation*}
    Then one has:
    \begin{equation}
      \sin^2\angle(m_i,u_j)+\sin^2\angle(n_i,v_j)\leq{{\p R_j\p^2+\p S_j\p^2}\over{\de_{i,j}^2}}.
    \end{equation}
    Note: the above statement is a special case of a more generic form of Wedin's Theorem given in \cite[p.260]{Stewart1990}. It is stated there for a generic complex case, whereas the above form is obtained simply by restricting its proof to real case and one-dimensional subspaces spanned by singular vectors $m_i,n_i,u_j,v_j.$
  \end{fact}

  \begin{fact}[\cite{Dopico2000}] \label{Dopico}
    Let $A$ and $\widetilde{A}=A+E$ for some $E\in\nm$ be two $n\times m$ matrices such that $n\geq m$, both with distinct singular values. For $i,j\in\{1,\dots,m\}$, let us pick $\alpha_i=i$-th singular value of $A$ along with $m_i$ and $n_i$: left and right singular vectors corresponding to $\alpha_i$, and similarly $\tilde{\alpha}_j=j$-th singular value of $\tilde{A}$ along with $u_j$ and $v_j$: left and right singular vectors corresponding to $\tilde{\alpha}_j.$ Let us define:
    \begin{equation}
      \zeta_{i,j}=\min\{\de_{i,j},\alpha_i+\tilde{\alpha}_j\},
    \end{equation}
    where $\de_{i,j}$ is defined as in (\ref{Wedin_delta}). Moreover, let us set:
    \begin{equation*}
      R_j=(A-\widetilde{A})v_j
    \end{equation*}
    and 
    \begin{equation*}
      S_j=(A^t-\widetilde{A}^t)u_j.
    \end{equation*}
    If $\zeta_{i,j}>0$, then one has:
    \begin{equation} \label{hasDopico}
      \min_{x\in\{-1,1\}}\left\{\p xm_i-u_j\p^2+\p xn_i-v_j\p^2\right\}\leq 2{{\p R_j\p^2+\p S_j\p^2}\over{\zeta_{i,j}^2}}.
    \end{equation}
    Note: the above statement is a special case of a more generic form of Theorem 2.1. in \cite{Dopico2000}. It is stated there for a generic complex case, whereas the above form is obtained simply by restricting its proof to real case and one-dimensional subspaces spanned by singular vectors $m_i,n_i,u_j,v_j.$
  \end{fact}

  \section*{Acknowledgment} 
  This study as a part of ongoing research project was supported by a grant from the Polish National Science Centre (UMO- 2013/08/W/HS6/00333) “NeuroPerKog: development of phonematic hearing and working memory in infants and children”, awarded to Prof. W Duch. The authors are grateful to anonymous reviewers for their constructive comments which surely promoted the readability of the revised manuscript. The first author would also like to thank Dr. J Dreszer for her continuous support and encouragement. 
  \section*{References} 

  \bibliographystyle{model1-num-names}
  \bibliography{IEEEabrv,references}







\end{document}
